\documentclass[12pt] {amsart}

\usepackage{amssymb,amsmath, amsthm,latexsym}
\usepackage{a4wide}

\parskip=\medskipamount
\font\k=cmr7
\font\rm=cmr12
\newcommand{\length}{\operatorname{-length}}
\newcommand{\slimit}{\operatorname{s-lim}}
\newcommand {\C}{\mathbb C}
  \newcommand {\N}{\mathbb N}
  \newcommand {\R}{\mathbb R}

\newcommand {\cA}{{\mathcal A}}
 \newcommand {\cB}{{\mathcal B}}
\renewcommand {\H}{{\mathcal H}}
  \newcommand {\M}{{\mathcal M}}
  
  \newcommand {\cF}{{\mathcal F}}
\newcommand {\Ell}{{\mathcal Ell}}
\newcommand {\sD}{\sqrt{\Delta}}

\newcommand {\p}{\parallel}
\newcommand{\loc}{\mbox{\k loc}}

\renewcommand {\Re}{\operatorname{Re}}
\newcommand {\Id}{\operatorname{Id}}
\newcommand {\ac}{{\mathrm ac}}
\newcommand {\supp}{\operatorname{supp}}
\newcommand {\Hom}{\operatorname{Hom}}
\newcommand {\vol}{\operatorname{vol}}
\newcommand {\ov}{\overline}
\newcommand {\Tr}{\qopname\relax o{Tr}}

\theoremstyle{plain}
\newtheorem{lem}{Lemma}[section]
\newtheorem{theo}[lem] {Theorem}
\newtheorem{cor}[lem]{Corollary}

\newtheorem{prop}[lem]{Proposition}

\theoremstyle{definition}
\newtheorem{defn}[lem]{Definition}

\theoremstyle{remark}
\newtheorem{rem}{\bf Remark} 
\begin{document}
\date{\today}
\title[Scattering theory]{Scattering theory for the
  Laplacian  on manifolds with bounded curvature} 
\author{Werner M\"{u}ller}
\author{Gorm Salomonsen.}
\address{Universit\"at Bonn\\
Mathematisches Institut\\
Beringstrasse 1\\
D -- 53115 Bonn, Germany}
\address{Cryptomathic A/S\\
Jaegergaardsgade 118\\
DK-8000 Aarhus C\\
Denmark}
\email{mueller@math.uni-bonn.de}
\email{gorm.salomonsen@cryptomathic.com}


\begin{abstract}
In this paper we study the behaviour of the continuous spectrum of the 
Laplacian on a complete Riemannian manifold of bounded curvature under 
perturbations of the metric. The perturbations that we consider are such 
that its covariant derivatives  up to some order decay with some rate in
the geodesic distance from a fixed point. Especially
we impose no conditions on the injectivity radius. One of the main results 
are conditions on the rate of decay, depending on geometric properties of
the underlying manifold,  that guarantee the existence and completeness of
the wave operators.

\end{abstract}
\maketitle

\setcounter{section}{-1}
\section{Introduction.}
\setcounter{equation}{0}

Let $(M,g)$ be a complete Riemannian manifold and let $\Delta_g$ be the 
Laplacian on functions attached to $g$. Then $\Delta_g$ is an essentially
self-adjoint operator in $L^2(M)$ \cite{Cn}. If $M$ is non-compact, then 
$\Delta_g$ may have a nonempty continuous spectrum. The purpose of this 
paper is to study the behavior of the continuous spectrum of $\Delta_g$ 
with respect to perturbations of the metric $g$. If $h$ is a compactly
supported perturbation of $g$, 
it is well known that the wave operators 
$$W_{\pm}(\Delta_g,\Delta_h):=\slimit_{t\to\pm\infty}
e^{it\Delta_g}Je^{-it\Delta_h}P_{\mathrm{ac}}(\Delta_h)$$ 
exist.
Therefore, the absolute continuous parts $\Delta_{g,\mathrm{ac}}$ and
$\Delta_{h,\mathrm{ac}}$ of $\Delta_g$ and $\Delta_h$, respectively, are
unitarily equivalent. 
    Our goal is to study non-compactly supported perturbations of the metric.
To this end  we introduce a certain class of functions, called functions of
moderate decay, which  describe the rate of decay of the perturbation of a 
given metric.  Let $\beta\colon [1,\infty)\to \R^+$ be a function of moderate 
decay (see Definition \ref{d1.3}). Then two complete metrics $g$ and $h$
are said to be equivalent up to order $k\in\N$, if there exist $C>0$ and 
$p\in M$ such that
$$|g-h|_g(x)+\sum^{k-1}_{j=0}\left|(\nabla^g)^j(\nabla^g-\nabla^h)
\right|_g(x)\le C\beta(1+d_g(x,p)),\quad x\in M,$$
where $d_g(x,p)$ is the geodesic distance of $x$ and $p$ with respect to $g$, 
and $\nabla^g$ (resp. $\nabla^h$) the Levi-Civita connection with respect to 
$g$ (resp. $h$). This condition turns out to be an equivalence relation in the
set of complete metrics on $M$. We denote this equivalence relation by
$g\sim_\beta^k h$. It implies, in particular, that the two metrics are 
quasi-isometric.

To develop scattering theory for the Laplacian we need to impose
additional assumptions on the metrics. 
In this paper we restrict attention to the
class of complete metrics with bounded sectional curvature.
In some cases we will also demand that higher derivatives of the curvature
 tensor are bounded. The assumption that the metric has bounded sectional
curvature allows us to control the behavior of the injectivity radius 
$\imath(x)$ sufficiently well.

One of the main results is the following theorem.

\begin{theo}\label{th0.1}
Assume $g $ and $h$ be complete metrics on $M$ with bounded curvature up to
 order $2$. Let $\beta $ be a function of moderate decay. Suppose that 
$g\sim_\beta^2 h$.
Assume  that there exist real numbers $a,b$ satisfying   
\begin{itemize} 
\item[i)] $b\geq 1$ and $a+b=2$,
\item[ii)] $\beta^\frac{b}{3} \in L^1(M)$,
\item[iii)] $\beta^\frac{a}{3}\tilde{\imath}^{-\frac{n(n+2)}{2}} 
\in L^\infty(M)$.
\end{itemize}
Then $e^{-t\Delta_g}-e^{-t\Delta_h}$ is a trace class operator. 
\end{theo}

Here $\tilde\imath(x)$ is the modified injectivity radius, defined by 
(\ref{2.1}), which is bounded from above by a constant which depends on the
bound of the sectional curvature. Moreover, $e^{-t\Delta_h}$ is regarded as 
bounded operator in $L^2(M,g)$. This is possible, because $g$ and $h$ are 
quasi-isometric.

By the  invariance principle for wave operators \cite{Ka}, Theorem \ref{th0.1}
implies  that the
wave operators $W_\pm(\Delta_g,\Delta_h)$ exist and are complete (see Theorem 
\ref{th7.1}). Under additional assumptions on $(M,g)$, the conditions on
$\beta$ can be relaxed. This is, for example, the case for manifolds with 
cusps and manifolds with cylindrical ends. In either case, the method of Enss
can be used to prove the existence and completeness of the wave operators.

We also study the analytic continuation of the resolvent. This result can be
used, for example, to construct generalized eigenfunctions as in \cite{Mu2}.
 
The structure of the paper is as follows. In section 1 we introduce our class
of functions of moderate decay and study some of its elementary properties.
Then we set up the equivalence relation mentioned above and prove some facts 
about equivalent metrics. In section 2 we study the behavior of the injectivity
radius on manifolds with bounded sectional curvature. 
 
Then we introduce and study weighted Sobolev spaces in section 3. In section 4
we show that certain  functions of the Laplacian including the heat kernel and
the resolvent extend to bounded operators 
in weighted $L^2$-spaces. Section 5 deals with the comparison of weighted
Sobolev spaces with respect to equivalent metrics. Then we prove Theorem 
\ref{th0.1} in section 6. In section 7 we deal with the existence and 
completeness of wave operators. First we prove a general result which is
based on Theorem \ref{th0.1}. Then we consider the case of a manifold with
cusps and use the method of Enss to establish the existence and
completeness of the wave operators under weaker assumptions on $\beta$.
The final section 8 deals with the analytic 
continuation of the resolvent, regarded as operator in weighted $L^2$-spaces.

\section{Equivalence of Riemannian Metrics.}
\setcounter{equation}{0}
\sectionmark{EQUIVALENCE OF METRICS}

Let $M$ be an open, connected $C^{\infty}$--manifold of dimension $n$ and let
$\M=\M(M)$ be the space of all complete Riemannian metrics on $M$. 
Eichhorn \cite{Ei1} has shown that $\M$ can be endowed with a canonical
topology given by a metrizable uniform structure. We briefly recall its
definition. 

For a given Riemannian metric $g$ on $M,$ denote by $\nabla^g$ the Levi--Civita connection of $g$ and by $|\cdot|_g$ the norm
induced by $g$ in the fibers of $\oplus_{p,q\ge0}(TM^{\otimes p}\otimes
T^*M^{\otimes q}).$ Let $h$ be any other Riemannian metric on $M.$ For $k\ge0$
set

\begin{align}\label{1.1}
^k|g-h|_g(x) & =|g-h|_g(x)+\sum^{k-1}_{j=0}\left|(\nabla^g)^j(\nabla^g-\nabla^h)\right|_g(x),\quad x\in M. \\
\intertext{and}
^k|\hskip-.7pt| g-h\parallel_g & =\sup_{x\in M}{^k|g-h|_g(x)}.\end{align}
Recall that two metrics $g,h$ are said to be quasi-isometric if 
there exist $C_1,C_2>0$ such that

\begin{equation}\label{1.2}
C_1g(x)\le h(x)\le C_2g(x),\quad \text{for all }x\in M, 
\end{equation}
in the sense of positive definite forms. 
We shall write $g\sim h$ for quasi--isometric metrics $g$ and $h.$ If $g$ and
$h$ are quasi--isometric, then (\ref{1.2}) implies that for all $p,q\ge0,$ there exist
$A_{p,q},B_{p,q}>0$ such that for every tensor field $T$ on $M$ of bidegree
$(p,q),$ we have

\begin{equation}\label{1.3}
A_{p,q}|T|_g(x)\le|T|_h(x)\le B_{p,q}|T|_g(x),\quad x\in M.
\end{equation}
Put $\nabla:=\nabla^g$ and $\nabla^\prime:=\nabla^h$. Let $\nabla^{p,q}$ and
${\nabla^\prime}^{p,q}$ be the canonical extension of $\nabla$ and $\nabla^\prime$, respectively, to the tensor bundle
$T^{p,q}(M)$. Then for all $p,q\in \N$ there exists $C_{p,q}>0$ such that
\begin{equation}\label{1.3a}
|\nabla^{p,q}-{\nabla^\prime}^{p,q}|_g(x)\le C_{p,q}|\nabla-\nabla^\prime|_g(x),
\quad x\in M.
\end{equation}

\noindent For $k\ge1$ and $\delta>0,$ set
$$V_\delta=\{(g,g')\in\M\times\M\mid g\sim
g'\;\hbox{and}\;^k|\hskip-.7pt| g-g'\parallel_g<\delta\}.$$
It is proved in \cite{Ei1}, Proposition 2.1, that $\{V_\delta\}_{\delta>0}$ is a
basis for a metrizable uniform structure on $\M.$ 

\begin{lem} \label{l1.1a}  Let
$g,h\in\M$. Assume that there exists a compact subset $K\subset M$ and
$0<\delta<1$ such that
$|g-h|_g(x)\le\delta$ for all $x\in M\setminus K$. Then 
$g$ and $h$ are quasi-isometric.
\end{lem}
\begin{proof} Let $x\in M\setminus K$. Choose geodesic coordinates w.r.t. $g$,
centered at $x$. Then $g_{ij}(x)=\delta_{ij}$. Let $H=(h_{ij}(x))$ be the
matrix 
representing $h(x)$ in these coordinates. Denote by $\parallel\cdot\parallel$
the supremum norm of  linear maps in $\R^n$. Then by assumption, we have
$\parallel H-\text{Id}\parallel\le\delta<1$. Hence the Neumann series for
$H^{-1}=(\text{Id}-(\text{Id}-H))^{-1}$ converges in norm which implies that 
$\parallel H^{-1}\parallel\le1/(1-\delta)$. Thus for all $\xi\in\R^n$, we get
$$(1-\delta)\parallel\xi\parallel^2\le\left(\parallel
  H^{-1}\parallel\right)^{-1} 
\parallel\xi\parallel^2\le\langle H\xi,\xi\rangle\le\parallel H\parallel\;
\parallel\xi\parallel^2\le(1+\delta)\parallel\xi\parallel^2.$$
This implies that
$$(1-\delta)g(x)\le h(x)\le(1+\delta)g(x),\quad\text{for all } x\in M\setminus
K.$$ 
Since $K$ is compact, it follows that $g$ and $h$ are quasi-isometric.
\end{proof}

We need two results from the
proof of Proposition 2.1 in \cite{Ei1} which we state as lemmas. For the 
convenience
of the reader we repeat the proofs.
\begin{lem}\label{lem1.1} 
Let $g,h\in\M$ be quasi--isometric.
For every $k\ge0,$
there exists a polynomial $P_k(X_1,...,X_k)$,
depending on the quasi-isometry constants,  
with nonnegative coefficients and vanishing constant term, such that
$$^k|g-h|_h(x)\le P_k\bigl(|g-h|_g(x),|\nabla^g-\nabla^h|_g(x),
...,|(\nabla^g)^{k-1}(\nabla^g-\nabla^h)|_g(x)\bigr),\quad x\in M.$$
\end{lem}
\begin{proof}
From (\ref{1.3}) follows that
\begin{equation}\label{1.4}
|g-h|_h(x)  \le C_3|g-h|_g(x)
\end{equation}
and
\begin{equation}\label{1.4a}
|\nabla^g-\nabla^h|_h(x)  \le C_4|\nabla^g-\nabla^h|_g(x),\quad x\in M.
\end{equation}
This takes care of the first two terms in (\ref{1.1}) and settles the
question for
$k=0,1$. Now we shall proceed by induction. Let $k\ge2$ and suppose that the
lemma holds for $l\le k-1.$ For each $p\ge0,$
we have
\begin{equation}\label{1.5}
(\nabla^h)^p(\nabla^h-\nabla^g)=\nabla^g(\nabla^h)^{p-1}(\nabla^h-
\nabla^g)+(\nabla^h-\nabla^g)(\nabla^h)^{p-1}(\nabla^h-\nabla^g).
\end{equation}
Let $p\le k$. Using (\ref{1.4a}), (\ref{1.3a}) and the induction 
hypothesis, we can estimate the pointwise $h$-norm of the second term on the
 right hand side of (\ref{1.5}) in the desired way. To deal with the first
term, we use the formula
\begin{equation}\label{1.6}
\begin{split}
(\nabla^g)^p(\nabla^h)^l(\nabla^h-\nabla^g) =&
(\nabla^g)^{p+1}(\nabla^h)^{l-1}
(\nabla^h-\nabla^g)\\
&+(\nabla^g)^p(\nabla^h-\nabla^g)(\nabla^h)^{l-1}(\nabla^h-\nabla^g).
\end{split}
\end{equation}

\noindent Applying the Leibniz rule, we get

\begin{equation*}
\begin{split}
(\nabla^g)^p(\nabla^h-\nabla^g) (\nabla^h)^{l-1}(\nabla^h&-\nabla^g)\\
&=\sum^p_{i=0} \binom {p}{i} \bigl((\nabla^g)^i(\nabla^h-\nabla^g)\bigr)
\bigl((\nabla^g)^{p-i}
(\nabla^h)^{l-1}(\nabla^h-\nabla^g)\bigr).
\end{split}
\end{equation*}

\noindent Inserting (\ref{1.5}) and iterating these formulas reduces 
everything 
to the
induction hypothesis. 
\end{proof}

\begin{lem}\label{lem1.2}
Let $g_i\in\M,\; i=1,2,3,$ and suppose that $g_1
\sim g_2\sim g_3.$ For every $k\ge0,$ there exists a polynomial $Q_k$,
depending on the quasi-isometry constants, in the variables 
$^i|g_1-g_2|_{g_1}(x)$ and
$^j|g_2-g_3|_{g_2}(x),\quad i,j=0,\ldots,k,$ with nonnegative 
coefficients and vanishing constant term, such that
$$^k|g_1-g_3|_{g_1}(x)\le Q_k(^i|g_1-g_2|_{g_1}(x),\;
^j|g_2-g_3|_{g_2}(x)),\quad x\in M.$$
If there exists $\delta<1$ such that $\|g_1-g_2\|_{g_1}\le\delta$ and $\|g_2-g_3\|_{g_2}\le\delta$, 
the dependence on the quasi-isometry constants can be removed.
\end{lem}

\begin{proof}
Since $g_1\sim g_2,$ it follows from (\ref{1.3}) that
$$|g_1-g_3|_{g_1}(x)\le|g_1-g_2|_{g_1}(x)+C_1|g_2-g_3|_{g_2}(x).$$
Set $\nabla_i=\nabla^{g_i},\; i=1,2,3.$ By the same argument, we get
$$|\nabla_1-\nabla_3|_{g_1}(x)\le|\nabla_1-\nabla_2|_{g_1}(x)+
C_2|\nabla_2-\nabla_3|_{g_2}(x).$$
Thus, the lemma holds for $k=0,1,$ and we can use induction to prove the
lemma. First observe that for $p\ge0,$
$$\nabla^p_1(\nabla_1-\nabla_3)=\nabla^p_1(\nabla_1-\nabla_2)+
\nabla^p_1(\nabla_2-\nabla_3).$$
The pointwise $g_1$--norm of the first term on the right hand side gives
already what we want. The second term can be written as
$$\nabla^p_1(\nabla_2-\nabla_3)=(\nabla_1-\nabla_2)\nabla^{p-1}_1
(\nabla_2-\nabla_3)+\nabla_2\nabla_1^{p-1}(\nabla_2-\nabla_3).$$
Iteration of this formula and application of the Leibniz rule reduces again
everything to the induction hypothesis. The last statement again follows 
from Lemma \ref{l1.1a}.
\end{proof}

To set up our equivalence relation in $\M$, we introduce an appropriate class
of functions. 

\begin{defn}\label{d1.3} 
Let $\beta: [1,\infty)\rightarrow\R$ be a positive, continuous, non-increasing
function. 
 Then $\beta$ is called a function of {\it moderate decay}, if it satisfies the
following conditions
\begin{equation}\label{1.7}
\begin{split}
&1)\;\; \sup_{x\in[1,\infty)}
x\beta(x)<\infty ;\\
 &2)\;\;\exists \;C_\beta>0:
\beta(x+y)\ge C_\beta\;\beta(x)\beta(y), \quad x,y\ge 1.
\end{split}
\end{equation}
Furthermore, $\beta$ is called of {\it  sub-exponential decay} if 
for any $c>0$, $e^{cx}\beta(x)\to\infty$ as $x\to\infty$.
\end{defn}

\begin{rem}
The class of functions which are of moderate or sub-exponential decay
are closed under multiplication, and also under raising to positive powers.
The function $e^{-tx}$, $t\ge0$, is of moderate decay and the functions
$x^{-1}$ and $\exp(-x^\alpha)$, $0<\alpha<1$, are of sub-exponential decay. 
Thus the
class of functions introduced in Definition \ref{d1.3} is not empty. 
\end{rem}

Next we  establish some elementary properties of $\beta.$

\begin{lem}
Let $\beta$ be of moderate decay. Then there exist constants $C>0$ and
$c\ge0$ such that
\begin{equation}\label{1.7a}
\beta(x)\ge C e^{-cx},\quad x\in[1,\infty).
\end{equation}
\end{lem}
\begin{proof}
Given $x\in [1,\infty)$, write $x$ as $x=y+n$, where $y\in[1,2)$ and
$n\in\N$. Applying  condition 2) of (\ref{1.7}) repeatedly, we get 
\begin{equation}\label{1.7b}
\beta(x)\ge \beta(y)(C_\beta\beta(1))^n.
\end{equation}
By assumption, $\beta$ is continuous. Hence there exists $C>0$ such that 
$\beta(y)\ge C$ for
$y\in [1,2]$. 
Since $\beta$ is non-increasing, it follows that $C_\beta\beta(1)\le 1$. Thus
there exists $c\ge 0$ such that $C_\beta\beta(1)=e^{-c}$. Together with
(\ref{1.7b}) the claim follows.
\end{proof}
Thus for a function $\beta$ of moderate decay there exist constants $c,C_1,C_2>0$
such that
$$C_1e^{-cx}\le\beta(x)\le C_2 x^{-1},\quad x\ge 1.$$
\begin{lem}\label{l1.4}
Let $\beta$ be a function  of moderate decay. Then for all $x,y,q\in M$, we
 have
\begin{equation}\label{1.8}
C_\beta\;\beta\bigl(1+d(x,y)\bigr)\le
{\frac{\beta\bigl(1+d(x,q)\bigr)}{\beta\bigl(1+d(y,q)\bigr)}}\le
\frac{1}{C_\beta\beta\bigl(1+d(x,y)\bigr)}. 
\end{equation} 
Moreover, for every $q'\in M$ there exists a constant $C>0$, depending only
 on $q$ and $q'$, such that
\begin{equation}\label{1.9} 
C^{-1}\;\beta\bigl(1+d(x,q')\bigr)\le \beta\bigl(1+d(x,q)\bigr)\le C\;\beta\bigl(1+d(x,q')\bigr).
\end{equation}
\end{lem}

\begin{proof}
Since $\beta$ is non-increasing, it follows from (\ref{1.7}) that
\begin{equation*}
\begin{split} \frac{\beta\bigl(1+d(x,q)\bigr)}{\beta\bigl(1+d(y,q)\bigr)}&\le
\frac{\beta\bigl(1+d(x,q)\bigr)}{\beta\bigl(1+d(x,q)+1+d(x,y)\bigr)}\\
&\le \frac{\beta\bigl(1+d(x,q)\bigr)}{C_\beta\;\beta\bigl(1+d(x,q)\bigr)
\beta\bigl(1+d(x,y)\bigr)}\\
&=\frac{1}{C_\beta\beta\bigl(1+d(x,y)\bigr)}.
\end{split}
\end{equation*}
Switching the roles of $x$ and $y$, we obtain the other inequality in
(\ref{1.8}). Furthermore, switching the roles of $x$ and $q$  and
putting  $y=q'$ in (\ref{1.8}) gives (\ref{1.9}).
\end{proof}

\begin{lem}\label{l1.7} Let $\beta$ be a function of moderate decay.
Let $g,h\in\M$, $q\in M$, and suppose that 
\begin{equation}\label{1.10}
|g-h|_g(x)\le \beta(1+d_g(x,q)),\quad x\in M. 
\end{equation}
Then $g$ and $h$ are quasi-isometric and there exist
constants $C_1,C_2>0$ such that
\begin{equation}\label{1.11a}
C_1d_g(x,y)\le d_h(x,y)\le C_2d_g(x,y),\quad x,y\in M,
\end{equation}
and
\begin{equation}\label{1.11}
C_1\beta(1+d_g(x,q))\le\beta(1+d_h(x,q))\le C_2\beta(1+d_g(x,q)),\quad x\in M.
\end{equation}
\end{lem}

\begin{proof}
Let $0<\delta<1$. 
From condition 1) of (\ref{1.7}) follows that there  exists
$r_0$ such that $\beta(1+r)\le\delta$ for $r\ge r_0.$  
Thus by Lemma \ref{l1.1a}, $g$ and $h$ are quasi-isometric and this implies
(\ref{1.11a}).  
To prove the second part, we first note that it follows from the proof of 
Lemma \ref{l1.1a} that 
$$d_h(x,q)\le (1+\beta(1+d_g(x,q))d_g(x,q),\quad d_g(x,q)\ge r_0.$$
Moreover, by condition 1) of (\ref{1.7}) there exists $C>0$ such that
$$\beta(1+d_g(x,q))d_g(x,q)\le C,\quad x\in M.$$
Then
using (\ref{1.7}), (\ref{1.11a}) and the assumption that $\beta$ is
 non-increasing, we get
\begin{equation*}
\begin{split}
\beta(1+d_h(x,q))\ge\beta(1+(1+\beta(1+d_g(x,q))d_g(x,q))
\ge C_\beta\beta(C)\beta(1+d_g(x,q)).
\end{split}
\end{equation*}
Switching the roles of $g$ and $h,$ we obtain the other inequality.
\end{proof}

\noindent Let $k\ge0$, and consider the following relation for metrics 
$g,h\in\M\colon$

\begin{equation}\label{1.12}\begin{split}
\hbox{\it There exist } &q\in M\hbox{\,and\,\,}C>0\;\hbox{\it such that for
  all\,\,} x\in M \;{\it we\; have}\\
&^k|g-h|_g(x)\le C\beta(1+d_g(x,q)).
\end{split}
\end{equation}

\begin{prop}\label{prop1.6} 
The relation (\ref{1.12}) defines an equivalence relation in
$\M.$
\end{prop}

\begin{proof}
Let $g,h\in\M$ and suppose that (\ref{1.12}) holds. 
Then by Lemma \ref{l1.7}, $g,h$ are
quasi--isometric. Then Lemma \ref{lem1.1} 
combined with (\ref{1.11}) implies that  
$$^k|g-h|_h(x)\le C_3\beta(1+d_g(x,y))\le C_4 \beta(1+d_h(x,q)).$$
Thus the relation (\ref{1.12})
 is symmetric. The
transitivity follows 
from Lemma \ref{lem1.2} and (\ref{1.11}). By Lemma~\ref{l1.4}, the relation is
independent of $q$.
\end{proof}

\noindent This justifies the following definition.

\begin{defn}
Let $\beta$ be a function of moderate decay. 
Two metrics $g,h\in\M$  are said to be
$\beta$--equivalent up to order $k$ if (\ref{1.12}) holds. In this case we
write
$g\sim^k_\beta h.$
\end{defn}

\smallskip
\noindent
{\bf Example 1.10.} 
Let $(M,g)$ be a complete Riemannian manifold which is 
Euclidean at
  infinity, that is, 
there 
exists a compact subset $K\subset M$ such that 
$(M\setminus K,g)$ is isometric to $\R^n\setminus B_r(0)$ for some $r>0$, 
 where $\R^n$ is equipped with its standard metric. Let $\beta(r)=r^{-a}$,
 $a>1$, and let $h$ be a complete Riemannian metric on $M$ such that
 $h\sim^k_\beta g$ for some $k\in\N$. Then $h|_{M\setminus K}$ may
 be regarded as metric on $\R^n\setminus B_r(0)$ and if $h_{ij}$ are the
 components of $h|_{M\setminus K}$ with respect to the standard
 coordinates $x_1,...,x_n\in\R^n$, then the condition $h\sim^k_\beta g$ is
equivalent to 
\begin{equation}\label{1.12a}
\Bigg|\frac{\partial^\alpha}{\partial x_1^{\alpha_1}\cdots\partial
  x_n^{\alpha_n}}\left(h_{ij}(x)-\delta_{ij}\right)\Bigg|\le C(1+\parallel
  x\parallel)^{-a}
\end{equation}
for all multindeces $\alpha$ with $|\alpha|\le k$ and all $x\in\R^n\setminus
  B_r(0)$. Such metrics are called asymptotically Euclidean.

\bigskip

\noindent
To simplify notation, we will write $\beta(x)$ in place of $\beta(1+d_g(x,q)).$ If
$g\sim^k_\beta h,$ it follows from Lemma \ref{1.5}, that we  may use both $d_g$ and $d_h$ in (\ref{1.12}).

Next we show that the $\beta$-equivalence can also be defined in a 
different manner. Namely we have the following proposition.

\begin{prop}\label{prop1.8} 
Let $k\geq0$ and let $g,h\in\M$. Then $g\sim^k_\beta h$ holds if and only if there exists $C_1>0$ such that
$$\sum_{i=0}^k\big|(\nabla^g)^i(g-h)\big|_g(x)\le C_1\beta(x),\quad x\in M.$$
\end{prop}

\begin{proof}
Let $g,h\in\M$. The lemma holds obviously for $k=0$. Let $k\ge1$. Recall that $\nabla^gg=0$ and $\nabla^hh=0$. Using this fact, we get

\begin{equation}\label{1.13}
\begin{split}
(\nabla^g)^k(g-h) &  =-(\nabla^g)^kh=-(\nabla^g)^{k-1}(\nabla^g-\nabla^h)h\\
  &  =-\sum^{k-1}_{i=0} \left( \begin{matrix} k-1 \\ i
\end{matrix}\right)
\left((\nabla^g)^i(\nabla^g-\nabla^h)\right)\bigl((\nabla^g)^{k-1-i}(h)\bigr).
\\
\end{split}
\end{equation}
Suppose that $^k|g-h|_g(x)\le C\beta(x)$, $x\in M$, for some constant
$C>0$. Then $|h|_g(x)\le C^\prime$ for some constant $C^\prime>0$. By
induction it 
follows from (\ref{1.3a}) and (\ref{1.13}) that  
\begin{equation}\label{1.14} 
\sum_{i=0}^k\big|(\nabla^g)^i(g-h)\big|_g(x)\le C_1\beta(x),\quad x\in M,
\end{equation}
for some constant $C_1>0$, depending on $C$ and $k$.

Now assume that (\ref{1.14}) holds. We observe that for any smooth vector
fields $X,Y,Z$, the following formula holds

\begin{equation}\label{1.15}
\begin{split}
h\bigl((\nabla^g_X-\nabla^h_X)Y,Z\bigr)= \frac{1}{2}\;\bigl\{\nabla^g_X(g-h)(Y,Z)\,&+\,\nabla^g_Y(g-h)(X,Z)\\
&\, -\nabla^g_Z(g-h)(X,Y)\bigr\}.\\
\end{split}
\end{equation}

\noindent From this formula we get
$$|\nabla^h-\nabla^g|_h\le C\,|\nabla^g(g-h)|_h.$$
Taking covariant derivatives of (\ref{1.15}) and using induction, we obtain 
$$^k|h-g|_h(x)\le C \sum_{i=0}^k\big|(\nabla^g)^i(g-h)\big|_h(x).$$
By (\ref{1.3}) and (\ref{1.14}), we get
$$^k|h-g|_h(x)\le C\,\beta(x),$$ 
and Lemma \ref{lem1.1} implies that
$$^k|g-h|_g(x)\le C_1\,\beta(x),\quad x\in M,$$
for some constant $C_1>0$.
\end{proof}

Thus, we may define $\beta$--equivalence  also  by requiring that (\ref{1.14}) holds for some constant $C_1$. It follows from the previous proposition that this gives rise to an equivalence relation.

Finally, we study the behavior of the curvature tensor and its covariant derivatives under $\beta$--equivalence. Given $g\in \M,$ denote by $R^g$ the curvature tensor of $g.$

\begin{lem}\label{lem1.9} 
Let $k\ge2$ and let $g,h\in\M$. Suppose that $g\sim
^k_\beta h.$ Then there exists $C_k>0$ such that
$$|(\nabla^g)^i(R^g-R^h)|_g(x)\le C_k\beta(x),\;x\in M,\;i=0,\ldots,k-2.$$
\end{lem}

\begin{proof}
Set $\nabla=\nabla^g$, $\nabla'=\nabla^h$. We define the exterior differential
$$d^\nabla: C^\infty(\Lambda^p(T^\ast M)\otimes TM)\rightarrow 
C^\infty(\Lambda^{p+1}(T^\ast M)\otimes TM)$$
associated with $\nabla$ by the following formula
\begin{equation*}
\begin{split}
(d^\nabla\alpha)(X_0,...,X_p)=&\sum_{i=0}^p(-1)^i\nabla_{X_i}\bigl(\alpha(X_0,...,\widehat
X_i,...,X_p)\bigr) \\ 
&-\,\sum_{i<j}(-1)^{i+j}\alpha([X_i,X_j],X_0,...,\widehat X_i,...,\widehat
X_j,...,X_p).\\ 
\end{split}
\end{equation*}
Then, regarded as operators $C^\infty(TM)\rightarrow C^\infty(\Lambda^2(T^\ast
M)\otimes TM)$, we have 
$$R^\nabla=d^\nabla\circ d^\nabla,$$
and a corresponding formula holds for $\nabla'$. Set $A=\nabla'-\nabla$ and 
let $X,Y$ be smooth vector fields on $M$. Then we have \cite[p. 25]{Be}
\begin{equation*}
\begin{split}
R^{\nabla'}(X,Y)\;-\;R^\nabla(X,Y)&=\nabla_X\bigl(A(Y)\bigr)-\nabla_Y\bigl(A(X)\bigr)-A([X,Y])\\
&\hskip60pt-A(X)\circ A(Y)+A(Y)\circ A(X)\\
&=(\nabla A)(X,Y)-(\nabla A)(Y,X)\\
&\hskip60pt -A(X)\circ A(Y)+A(Y)\circ A(X).\\
\end{split}
\end{equation*}
Differentiating this equality and using induction gives the desired result. 
\end{proof}

Recall that a Riemannian manifold $(M,g)$ is said to have bounded curvature of order $k$, if the covariant derivatives $\nabla^iR$, $0\le i\le k$, of the curvature tensor $R$ are uniformly bounded on $M$, i.e., there exists $C>0$ such that $|\nabla^iR|(x)\le C$, $x\in M$, $0\le i\le k$.

\begin{cor}\label{cor1.10}
Let $k\ge2$ and let $g,h\in\M$. Suppose that $g\sim^k_\beta h$. Then
\begin{itemize}
\item[1)] $(M,g)$ has bounded curvature of order $k-2$ if and only if $(M,h)$
  has bounded curvature of order $k-2$.  
\item[2)] The sectional curvature of $(M,g)$ is bounded from below (resp. from
  above) if and only if the sectional curvature of $(M,h)$ is bounded from
  below (resp. above). 
\item[3)] The Ricci curvature of $(M,g)$ is bounded from below (resp. from
  above) if and only if the Ricci curvature of $(M,h)$ is bounded from below
  (resp. above). 
\end{itemize}
\end{cor}

\section{Injectivity radius and bounded curvature.} \label{geometrysection}
\setcounter{equation}{0}
\sectionmark{GEOMETRY OF MANIFOLDS WITH BOUNDED CURVATURE.}

In this section we establish some properties  of the
 injectivity radius  on a  manifold with bounded sectional curvature. 
Let $(M,g)$ be a complete, $n$--dimensional Riemannian manifold 
with bounded sectional curvature, say $|K_M|\le K.$ Let $p\in M$.
Recall that the injectivity radius $i(p)$ at $p$ equals the minimal 
distance from $p$ to its cut locus $C(p)$ (see \cite{CE}, \cite{Kl}).
Also note that $i(p)$ is a continuous function of $p\in M$
\cite[Proposition 2.1.10]{Kl}.
\begin{prop} \label{prop2.1}
Let $h$ be another complete Riemannian metric on $M$ with bounded 
sectional
curvature $|K_M^h|\le K$ and assume that $g$ and $h$ are equivalent.
 Given $p\in M$, let $i_g(p)$ and $i_h(p)$ denote the injectivity radii at $p$ 
with respect to $g$ and $h$, respectively. Then there exist constants 
$c,c'>0$ such that
$$i_h(p)\ge \min\{ci_g(p),c'\},\quad p\in M. $$
\end{prop}

\begin{proof}
Since $g$ and $h$ are assumed to be equivalent, there exists $\varepsilon>0$ 
such that
$$e^{-\varepsilon}g\le h\le e^{\varepsilon}g.$$
Let $x\in M$ and suppose that  $i_h(x)<\min\left\{
e^{-2\epsilon}\pi/ (2\sqrt{K}),e^{-\epsilon}i_g(x)/2\right\}$. It follows from
\cite[Corollary 1.30]{CE} that distinct conjugate points 
along a geodesic (with respect to $h$) have distance $\ge\pi/\sqrt K$. 
Therefore, by  \cite[Lemma 5.6]{CE},  there exists a closed geodesic loop 
$\gamma^h$ at $x$ with respect to the metric $h$,  with 
$$h\length(\gamma^h)<\min\left\{e^{-2\epsilon}\pi/\sqrt{K},
 e^{-\epsilon} i_g(x)\right\}.$$ 
Hence, we have
$$g\length(\gamma^h)<\min\left\{e^{-\epsilon}\pi/\sqrt{K},
\;i_g(x)\right\}.$$
In particular, $g$-length$(\gamma^h)<\pi/\sqrt K$.
Let $r_{\max}$ be the maximal rank radius of $\exp_x$ with respect to $g$. Then
we obtain $g$--length$(\gamma^h)<\pi/\sqrt K\le r_{\max}(x).$ By \cite{BK},
Proposition 2.2.2,
there exists a unique $g$-geodesic loop
${\widetilde{\gamma}}\colon\left[0,1\right]\longrightarrow M$ at $x$ with
$g$--length$({\widetilde{\gamma}})< r_{\max}(x),$ which is obtained from
$\gamma^h$ by 
a length decreasing homotopy
$H\colon\left[0,1\right]\times\left[0,1\right]\longrightarrow M$
(cf. \cite{BK}, 2.1.2). Hence, we have 
$$g\length({\widetilde{\gamma}})\le
g\length(\gamma^h)<\min\left\{ e^{-\epsilon}\frac{\pi}{\sqrt K},\;
 i_g(x)\right\}.$$
 Since $h$-length$(H(\cdot,s))\le e^\epsilon g$-length$(\gamma^h)<2\pi/\sqrt K
 $ for 
$s\in\left[0,1\right],$ it follows from \cite[Lemma 2.6.4]{Kl}, that
$g$-length$({\widetilde{\gamma}})>0.$
Parameterize ${\widetilde{\gamma}}$ by $g$--arc length. Then either
${\widetilde{\gamma}}(t)$ or
$\widetilde{\gamma}(\hbox{length}(\widetilde{\gamma})-t)$ belongs to the cut
locus of $x$ for some $t\le\frac{1}{2}g$-length$(\widetilde\gamma).$ 
Therefore $i_g(x)<i_g(x),$ a contradiction.
\end{proof}

Let $\beta$ be a function of moderate decay. 
Suppose that $g\sim^0_{\beta}h$. Then by Lemma \ref{l1.7}, $g$ and $h$ are
quasi-isometric.  Therefore, if $h$ has bounded sectional curvature, then 
Proposition \ref{prop2.1} can be applied to $g,h$. For $x\in M$ set
\begin{equation}\label{2.1}
\tilde{\imath}(x):=\min\left\{ \frac{\pi}{12\sqrt K},\;i(x)\right\}.
\end{equation}

\noindent Then it follows, that under the assumptions of
Proposition~\ref{prop2.1}, there exists $c_2>0$ such that
\[ \tilde{\imath}_h(p)\geq c_2\tilde{\imath}_g(p),\quad p\in M. \]

\noindent Next recall the Bishop--G\"unther inequalities \cite[Theorem
3.17]{Gra}, \cite[Lemma 5.3]{Gro}, which give estimates of the volume of small
balls from above and below. 

\begin{lem}\label{lem2.3} 
For $r\le\tilde{\imath}(x_0),$
$${ \frac{2\pi^{n/2}}  {\Gamma\left( \frac{n}{2}\right)} }
\int^r_0\left( \frac {\sin t\sqrt K}  {\sqrt K}
\right)^{n-1}dt\le\hbox{\rm Vol}(B_r(x_0))
\le { \frac {2\pi^{n/2}} {\Gamma\left(\frac{n}{2}\right)}}
\int^r_0\left( \frac {\sinh t\sqrt K} {\sqrt K}\right)^{n-1}dt.$$ 
\end{lem}

\noindent We note that the inequality on the right hand side holds for all
$r\in\R_+.$ In 
particular

\begin{equation}\label{2.1a}
\hbox{Vol}(B_r(x_0))=O\left(e^{(n-1)\sqrt K r}\right)
\end{equation}

\noindent as $r\to\infty.$

It is also important to know the maximal possible decay of the
injectivity radius.

\begin{lem} \label{lem2.4}
There exists a constant $C>0$, depending only on $K,$ such that
\begin{equation}\label{2.2}
\tilde{\imath}(x)\ge C\;\tilde{\imath}(p)^ne^{-(n-1)\sqrt K\;d(x,p)}
\end{equation}
for all $x,p\in M.$
\end{lem}

\begin{proof}
Let $p\in M$ and fix $r,r_0,s,$ with $r_0+2s<\pi/
\sqrt K,\;r_0\le\pi/4\sqrt K.$ By \cite[Theorem 4.7]{CGT} we get

\begin{equation} \label{e2.2}
\tilde{\imath}(x) \ge \frac{r_0}{2}\cdot
\frac {1}
{1+
\bigl(V^K_{r_0+s}/\hbox{Vol}(B_r(p))\bigr)
\left(V^K_{d(x,p)+r} /V^K_s \right)},
\end{equation}

\noindent where $V^K_s$ denotes the volume of a ball of radius $s$ in the
$n$--dimensional hyperbolic space of curvature $-K$. Set
$r_0=s= \frac {\pi}{5\sqrt K},\; r=\tilde{\imath}(p)$ and apply Lemma
\ref{lem2.3} 
 to
estimate $\hbox{Vol}\left(B_{\tilde i(p)}(p)\right)$ from below. Then
(\ref{e2.2}) 
implies

\begin{equation*}
\begin{split}
  \tilde{\imath}(x) &   \ge C_1\;\tilde{\imath}(p)^ne^{-(n-1)\sqrt
    K(d(x,p)+\tilde{\imath}(p))}\\ 
&  \ge C\;\tilde{\imath}(p)^n e^{-(n-1)\sqrt Kd(x,p)}.\\
\end{split}
\end{equation*}
\end{proof}

\begin{cor}\label{cor2.5} Given $p\in M,$ there exists a constant 
$C=C(p)>0$
  such that 
\begin{equation*}
\tilde{\imath}(x)\ge C e^{-(n-1)\sqrt K\;d(x,p)},\quad x\in M.
\end{equation*} \qed
\end{cor}

\begin{lem} \label{injbound2} There exists a constant $C$, 
depending only on
  $K$, such that for each $x,y\in M$ we have the inequality 
\begin{equation} \tilde{\imath}(y) \geq
  C\tilde{\imath}(x)e^{-\frac{(n-1)\pi}{12}\frac{d(x,y)}{\tilde{\imath}(x)}}.
  \end{equation}   
\end{lem}
\begin{proof} Let $\lambda = \max\{1,\frac{\pi^2}{144Ki(x)^2}\}$. 
Then the
  injectivity radius $i_\lambda$ at $x$ with respect to $\lambda g$ 
is given
  by 
\[ i_\lambda(x) = \lambda^\frac{1}{2}i(x) 
= \left\{ \begin{array}{ll} i(x) &
  ,\text{if }i(x)>\frac{\pi}{12\sqrt{K}}; \\ \frac{\pi}{12\sqrt{K}} & ,\text{if
  } 
i(x) \leq
  \frac{\pi}{12\sqrt{K}}. \end{array} \right. \] 
Since $\lambda^{-1} \leq 1$, the sectional curvature  $K_M^{\lambda g}$
with respect to 
$\lambda g$ also satisfies $|K_M^{\lambda g}|\le K$.

Let $r=\frac{\pi}{\sqrt{K}}$, $r_0=s=\frac{r}{12}
=\frac{\pi}{12\sqrt{K}}$ and
set $d=\lambda^\frac{1}{2}d_g(x,y)$. Then $d$ is the distance 
between $x$ and
$y$ with respect to $\lambda g$. 

Let $V_s(y)$ be the volume of the geodesic ball of radius $s$ and  
center $y$
with respect to $\lambda g$ and let $V_s^K$ denote the volume of a 
ball of radius $s$ in the
$n$-dimensional simply connected space of constant curvature $-K$. Then by
\cite[Theorem 4.3]{CGT} we get
\begin{equation} i_\lambda(y) \geq
  \frac{r_0}{2}\frac{1}{1+\frac{V^K_{r_0+s}}{V_s(y)}} \geq
  \frac{r_0}{4}\frac{V_s(y)}{V^K_{r_0+s}}. 
\label{ilambdaest} \end{equation}
Now, \cite[Proposition 4.1, i)]{CGT} states that
\[ \frac{V_s(y)}{V_s^K} \geq \frac{V_{d+s}(y)}{V_{d+s}^K}. \]
Together with \hbox{(\ref{ilambdaest})} this gives
\[ i_\lambda(y) \geq 
\frac{r_0}{4}\frac{V_{d+s}(y)V_s^K}{V_{d+s}^KV^K_{r_0+s}}. \]
From the definition of $d$ it  follows that, with respect to the metric
$\lambda g$,  the ball of radius $d+s$ 
around $y$
contains the ball of radius $s$ around $x$. Hence $V_{d+s}(y)\ge V_s(x)$.
Since
$s=\frac{\pi}{12\sqrt{K}}=\tilde{\imath}_\lambda(x)$, it follows from Lemma
 \ref{lem2.3}
that there exists $c>0$ such that $V_s(x)\ge c$ for all $x\in M$. Hence, we get
\begin{equation*}
\begin{split} i_\lambda(y) &\geq
  \frac{r_0}{4}\frac{V_{s}(x)V_s^K}{V_{d+s}^KV^K_{r_0+s}} \geq C
  \frac{V_s^K}{V_{d+s}^K}\geq  
Ce^{-(n-1)\sqrt{K}d} \\
&\geq Ce^{-(n-1)\max\{ \sqrt{K},\frac{\pi}{12i(x)} \}d(x,y)} 
= Ce^{-\frac{(n-1)\pi}{12}\frac{d(x,y)}{\tilde{\imath}(x)}}, 
\end{split}
\end{equation*}
for some constant $C>0$.
Now the lemma  follows by dividing both sides of this inequality by 
$\lambda^{\frac{1}{2}}$.
\end{proof}

We can now establish the following basic result about the existence 
of uniformly locally finite coverings on manifolds with bounded 
curvature. 

\begin{theo} \label{coverthm} Assume that $M$ is non-compact. Let $h$ be a
  continuous real valued function on 
  $M$ such that 
\begin{itemize}
\item[i)] $\forall x: 0< h(x) \leq \tilde{\imath}(x)$.
\item[ii)] There exists constants $C_1,C_2>0$ such that
\[ h(x) \geq C_1h(x_0)e^{-C_2\frac{d(x,x_0)}{h(x_0)}} \]
for all $x,x_0\in M$.
\end{itemize}
Then for each $a\geq 1$, there exists a sequence 
$\{ x_i\}_{i=0}^\infty\subset M$
 and a constant $C_3 <\infty$, depending only on $K$, $a$, $C_1$ and $C_2$
 such that 
\begin{enumerate}
\item[1)] \begin{equation*} \bigcup_{i=0}^\infty B_{h(x_i)}(x_i) =
    M. \end{equation*} 
\item[2)]  $\forall i\in\N:\# \{j\mid B_{ah(x_i)}(x_i)\cap B_{ah(x_j)}(x_j) \not= \emptyset \} \leq C_3$.
\end{enumerate}  
\end{theo}
\begin{proof} Let $x_0\in M$. For $k\in\N$  define recursively 
$$m(k)=\min\{m\in\N\mid
B_m(x_0)\setminus\cup_{i<k}B_{h(x_i)}(x_i)\not=\emptyset\}$$ 
and pick $x_k\in B_{m(k)}\setminus\cup_{i<k}B_{h(x_i)}(x_i)$.
In this way we get a sequence $\{x_i\}_{i=0}^\infty$ of points of $M$. From
the construction it follows that this sequence  satisfies the following
condition:
\begin{equation}\label{2.5}
\forall i,j\in\N: \;d(x_i,x_j)\ge\min\{h(x_i),h(x_j)\}.
\end{equation}
Let $m\in\N$. Then by ii), there exists $c>0$ such that $h(x)\ge c$ for all 
$x\in B_m(x_0)$. Hence it follows from (\ref{2.5}) that $d(x_i,x_j)\ge c$ if 
$x_i,x_j\in B_m(x_0)$. Since $B_m(x_0)$ is compact, this implies that only
finitely many of the $x_i$'s, say $x_1,...,x_{r_m}$, are contained in
$B_m(x_0)$. Hence 
$$B_m(x_0)\subset\bigcup_{i=0}^{r_m}B_{h(x_i)}(x_i)$$
which implies that
$$M=\bigcup_{i=0}^\infty B_{h(x_i)}(x_i).$$
It remains to prove 2). Let $a\ge1$. For $j\in\N$ put
$B_j=B_{ah(x_j)}(x_j)$.  Let  $i\in\N$ be given and put
$$\Omega_i=\{x_j\mid B_i\cap B_j\not=\emptyset\}.$$
Since $h$ is bounded from above, $\Omega_i$ is contained in a compact subset
$Y$ of $M$. By ii) there exists $c>0$ such that $h(x)\ge c$ for all $x\in Y$.
Using (\ref{2.5}), it follows that $\Omega_i$ is a discrete subset of $Y$ and
hence, 
$\Omega_i$ is a finite set. Let $x_{j_1}\in\Omega_i$ be such that
$$h(x_{j_1})=\max\{h(x_j)\mid x_j\in\Omega_i\}.$$
Since $B_i\cap B_{j_1}\not=\emptyset$, it follows that $B_i\subset
B_{3ah(x_{j_1})}(x_{j_1})$ which in turn implies that
$$B_{\frac{h(x_j)}{2}}(x_j)\subset B_{(4a+1)h(x_{j_1})}(x_{j_1})$$
for all $x_j\in\Omega_i$. Therefore  by ii) we get 
\[ h(x_{j}) \geq C_1 h(x_{j_1})e^{-C_2\frac{d(x_{j_1},x_j)}{h(x_{j_1})}}
 \geq C_1 h(x_{j_1})e^{-4aC_2}. \]
Thus there exists $C_4>0$ such that
\begin{equation}\label{2.6}
 h(x_j) \geq C_4h(x_{j_1}) 
\end{equation}
for all $x_j\in\Omega_i$.  Obviously  $C_4\leq 1$. Hence by i), we obtain
\begin{equation} \frac{C_4h(x_{j_1})}{2}\le\frac{\tilde{\imath}(x_{j_1})}{2}
 \leq \frac{\pi}{24\sqrt{K}}. \label{hbd} \end{equation}
Moreover, by (\ref{2.5}) and (\ref{2.6}) we have $d(x_k,x_l)\ge C_4
h(x_{j_1})$ for all $x_k,x_l\in\Omega_i$. Therefore,  the  balls
$B_{\frac{C_4}{2}h(x_{j_1})}(x_j)$, $x_j\in\Omega_i$, are pairwise 
disjoint.
Using Lemma \ref{lem2.3}, we get
\begin{equation} \# \{ x_j\mid B_i\cap B_j\not=\emptyset \} \leq
  \frac{\int_0^{(4a+1)h(x_{j_1})} \left( \frac{\sinh t\sqrt{K}}{\sqrt{K}}
    \right)^{n-1} dt}{\int_0^{\frac{C_4 h(x_{j_1})}{2}} \left( \frac{\sin
        t\sqrt{K}}{\sqrt{K}} \right)^{n-1} dt}. \label{numberest}
\end{equation} 
There exist 
constants $c_1>0$ and $c_2>0$ such that 
\[ \sinh t\sqrt{K} \leq c_1 t, \quad 0\le
t\le\frac{(4a+1)\pi}{12C_4\sqrt{K}};\] 
\[ \sin t\sqrt{K} \geq c_2 t,\quad 0\le t\le\frac{\pi}{24\sqrt{K}}. \]
Hence by (\ref{hbd}), it follows that the right hand side of (\ref{numberest})
is bounded by  
$\frac{c_2}{c_1}\left( \frac{(8a+2)c_1}{C_4c_2}\right)^n.$
 This proves the lemma. \end{proof}

\noindent Finally we will define and estimate some global invariants of
$(M,g)$. 

\begin{defn} \label{kappadef} Let $s >0$. For $s>\varepsilon\geq 0$,
let $\kappa_\varepsilon(M,g;s)\in \N\cup\{\infty \}$ be the smallest number
such that there exists a sequence $\{ x_i \}_{i=1}^\infty$ such that 
$\{B_{s-\varepsilon}(x_i) \}_{i=1}^\infty$ is an open covering of $M$ and 
\begin{equation} \sup_{x\in M} \#\{i\in \N \mid x\in B_{3s+\varepsilon}(x_i)
  \} \leq \kappa_\varepsilon(M,g;s). \end{equation} 
Further, let $\kappa(M,g;s)= \kappa_0(M,g,s)$. Put $\kappa(M,g,0)=1$.
\end{defn}

\begin{lem} \label{lordjagged} $\kappa_\varepsilon(M,g;s)$ is finite for all
  $s>\varepsilon$. Moreover, there exist constants $C,c>0$, which depend only
on $K$, such that for
  $s>\frac{2\pi}{\sqrt{K}}+\varepsilon$, we have 
\[ \kappa_\varepsilon(M,g;s) \le Ce^{cs}. \]
\end{lem}
\begin{proof} We may proceed as  in the proof of  Theorem~\ref{coverthm}
 and  construct a sequence $\{ x_i\}_{i=1}^\infty\subset M$ such that 
$d(x_i,x_j)\geq s-\varepsilon$ for all $i,j\in \N$ and 
$\{B_{s-\varepsilon}(x_i)\}_{i=1}^\infty$ is a covering of $M$. Let $x\in M$.
 If $x \in B_{3s+\varepsilon}(x_i)$, it follows that 
$B_{\frac{s-\varepsilon}{2}}(x_i) \subset B_{5s}(x)$. 
Moreover,   $B_{\frac{s-\varepsilon}{2}}(x_i)\cap
B_{\frac{s-\varepsilon}{2}}(x_j)=\emptyset$ if $i\not=j$. Hence, we get
\begin{equation}\label{2.10} \#\{ i\mid x\in B_{3s+\varepsilon}(x_i)\} 
\leq \frac{\text{Vol}(B_{5s}(x))}{\min_i\text{Vol}(B_{\frac{s-\varepsilon}{2}}(x_i))}. 
\end{equation}
Next observe that for any $x_i$ with $d(x,x_i)\le 5s$ we have
$B_{5s}(x)\subset B_{10s}(x_i)$. Moreover, by Lemma
5.3 of \cite{Gro}, we 
have
\begin{equation*}
\frac{\text{Vol}(B_{10s}(x_i))}{\text{Vol}(B_{\frac{s-\varepsilon}{2}}(x_i))}
\le\frac{\int_0^{10s}\left(\sinh t\sqrt{K}\right)^{n-1}\;dt}
{\int_0^{\frac{s-\varepsilon}{2}}
\left(\sinh t\sqrt{K}\right)^{n-1}\;dt}.
\end{equation*}
Then combined with (\ref{2.10}) we obtain
$$\#\{ i\mid x\in B_{3s+\varepsilon}(x_i)\} \le \frac{\int_0^{10s}
\left(\sinh t\sqrt{K}\right)^{n-1}\;dt}{\int_0^{\frac{s-\varepsilon}{2}}
\left(\sinh t\sqrt{K}\right)^{n-1}\;dt}.$$
If $(s-\varepsilon)/2\ge\pi/\sqrt{K}$,  the right hand side can be
estimated by $Ce^{cs}$ for certain constants $C,c>0$ depending on $K$.
\end{proof}

\section{Weighted Sobolev Spaces} \label{weighted}
\setcounter{equation}{0}

In this section we introduce   weighted Sobolev spaces on manifolds with
 bounded curvature.

Let $(M,g)$ be a Riemannian manifold. Let $\nabla$ be the Levi-Civita
connection of $g$ and let $\Delta=d^*d$ be the 
Laplacian on functions with respect to $g$. Let $\xi$ be a positive,
measurable function on $M$, which is finite  a.e. Given $m\in N_0$, and 
$p\in N$, we define
the weighted  $L^p$-space $L^p_\xi(M,TM^{\otimes m})$ by
$$L^p_\xi(M,TM^{\otimes m})=\bigl\{\varphi\in L^p_{\loc}(M,TM^{\otimes m})\mid 
\xi^{1/p}\varphi\in L^p(M,TM^{\otimes m})\big\}.$$
Then for $k\in\N$ we define the weighted Sobolev space $W^{p,k}_\xi(M)$ by
\begin{equation}\label{3.1}
W^{p,k}_\xi(M)  = \Bigl\{ f\in  L^p_\xi(M)\mid \nabla^mf\in 
L^p_\xi(M,TM^{\otimes m})\;\text{for all }m=1,...,k \Bigr\}, 
\end{equation}
where $\nabla$ is applied iteratively in the distributional sense
and the  norm of $f\in W^{p,k}_\xi(M)$ is given by
\begin{equation}\label{3.1a}
\parallel f\parallel_{W^{p,k}_\xi}
=\left(\sum_{i=0}^k\int_M|\nabla^if(x)|^p_g\xi(x)\;dv_g(x)\right)^{1/p}.
\end{equation}
Then $W^{p,k}_\xi(M)$ is a Banach space. In this paper we will only 
consider the case $p=2$. To simplify notation we shall write $W^k_\xi(M)$
in place of $W^{2,k}_\xi(M)$. The closure of $C^\infty_0(M)$ in
$W^k_\xi(M)$ will be denoted by $W^k_{0,\xi}(M)$. 
We shall write  $W^{k}(M)$ for $W^{k}_1(M)$ and $W^k_0(M)$ for 
$W^k_{0,1}(M)$. Since 0 is not a weight, this cannot lead to any confusion.
Note that $W^k_\xi(M)$ and $W^k_{0,\xi}(M)$ are Hilbert spaces.
The weighted Sobolev space 
$H^l_\xi(M)$ is defined  for even integers $l$. Let $k\in\N$. Then
\begin{equation}\label{3.2}
H^{2k}_\xi(M)=\Bigl\{f\in L^2_\xi(M)\mid \Delta^lf\in L^2_\xi(M)\;
\text{for all }\;l=1,...,k \Bigr\}.
\end{equation}
The norm is given by
$$\parallel f\parallel^2_{H^{2k}_\xi}=
\sum_{j=0}^k\int_M|\Delta^jf(x)|^2\xi(x)\;dv_g(x).$$
As an  equivalent norm one can use the norm defined by
\begin{equation}\label{3.2b}
\parallel f\parallel_{H^{2k}_\xi}=\parallel (\Delta+\Id)^kf\parallel_{L^2_\xi}
.
\end{equation}
The closure of $C^\infty_0(M)$ in $H^{2k}_\xi(M)$ will be denoted by 
$H^{2k}_{0,\xi}(M)$. 
If $\xi\equiv1$, the Sobolev space  $H^{2k}_\xi(M)$ will be denoted 
by $H^{2k}(M)$ and $H^{2k}_{0,\xi}$ by $H^{2k}_0(M)$. 
Note that the Laplacian $\Delta$ induces a bounded operator 
\begin{equation}\label{3.2a}
\Delta_\xi\colon H^2_\xi(M)\to L^2_\xi(M)
\end{equation}
which is defined in  the obvious way.

Next we establish some elementary properties of weighted Sobolev spaces.

\begin{lem}\label{l3.1}
Assume that $\xi$ is continuous. Let $p,k\in\N$. Then
$C^\infty(M)\cap W^{p,k}_\xi(M)$ is dense in $W^{p,k}_\xi(M)$ and
 $C^\infty(M)\cap H^{2k}_\xi(M)$ is dense in $H^{2k}_\xi(M)$.
\end{lem}
\begin{proof} We proceed as in the proof of Theorem 1 in \cite[1.1.5]{Ma}.
Let $\{U_i\colon i\in I\}$ be a locally finite covering of $M$ such that
for each $i\in I$ there exists an open subset $V_i$ with $\ov U_i\subset V_i$
and $V_i$ is diffeomorphic to the unit ball in $\R^n$. 
Let $\{\varphi_i\colon i\in I\}$
be an associated partition of unity. Let $u\in W^{p,k}_\xi(M)$ and let 
$\varepsilon\in (0,1/2)$. For each $i\in I$
let $u_i=\varphi_i u$. Then $u_i$ belongs to $W^{p,k}_\xi(M)$ with 
$\supp u_i\subset U_i$. Since 
$\xi$ is continuous, it follows that $u_i\in W^{p,k}(U_i)$ and $\supp u_i$
is contained in the interior of $U_i$. 
Hence there exists a mollification 
$g_i\in C^\infty_c(U_i)$ of $u_i$  such that 
$$\parallel g_i-u_i\parallel_{W^{p,k}}\le \frac{\varepsilon^i}
{\max_{x\in\ov U_i}\xi(x)}.$$ 
\cite[Section 5.3]{Ev}. Then $$\parallel g_i-u_i\parallel_{W^{p,k}_\xi}\le
 \varepsilon^i.$$
Clearly $g=\sum_i g_i$ belongs to $C^\infty(M)$. Let
$\omega\subset M$ be a relatively compact open subset. Then  we have 
$$u|_{\omega}=\sum_iu_i|_{\omega},$$
and the sum is finite. Hence
$$\parallel g-u\parallel_{W^{p,k}_\xi(\omega)}\le \sum_i\parallel g_i-u_i
\parallel_{W^{p,k}_\xi}\le \varepsilon(1-\varepsilon)^{-1}\le 2\varepsilon.$$
This implies that $\parallel u\parallel_{W^{p,k}_\xi(\omega)}\le \parallel u
\parallel_{W^{p,k}}+2\epsilon$ for all relatively compact open subsets 
$\omega\subset M$. Hence by the theorem of Beppo-Levi, we have
 $g\in C^\infty\cap W^{p,k}_\xi(M)$ and 
$$\parallel g-u\parallel_{W^{p,k}_\xi}\le 2\varepsilon.$$
The proof that $C^\infty(M)\cap H^{2k}_\xi(M)$ is dense in
$H^{2k}_\xi(M)$ is similar.
\end{proof}

Therefore we can use the following alternative definition of the Sobolev 
spaces. Let $C^\infty_k(M)$ denote the space of all $f\in C^\infty(M)$ such
that 
$|\nabla^jf|\in L^p_\xi(M)$ for $j=0,...,k$. Then $W^{p,k}_\xi(M)$ is the
completion of $C^\infty_k(M)$ with respect to the norm (\ref{3.1a}). Similarly
let  
$\tilde C^\infty_k(M)$ the space of all $f\in C^\infty(M)$ such that 
$(\Delta+\Id)^kf\in L^2_\xi(M)$. Then $H^{2k}_\xi(M)$ is the completion of 
$\widetilde C^\infty_k(M)$ with respect to the norm (\ref{3.2b}). This implies
that 
we can define $H^s_\xi(M)$ for all $s\in\R$. Let $(\Delta+\Id)^{s/2}$ be 
defined by the spectral theorem. Let $\widetilde C^\infty_s(M)$ be the space
of all $f\in C^\infty(M)$ such that $(\Delta+\Id)^{s/2}f\in L^2_\xi(M)$.
Let $H^s_\xi(M)$ be the completion of $\widetilde C^\infty_s(M)$ with respect
to the norm 
$$\parallel f\parallel_{H^s_\xi(M)}:=\parallel
(\Delta+\Id)^{s/2}f\parallel_{L^2_\xi}.$$

In general the Sobolev spaces  $W^k_\xi(M)$ and $W^k_{0,\xi}(M)$ (resp.
$H^{2k}_\xi(M)$ and $H^{2k}_{0,\xi}(M)$) will not coincide. If $(M,g)$ is
complete and $\xi\equiv 1$, the following is known \cite{Sa} .

\begin{lem}\label{l3.2} Assume that $(M,g)$ is complete. Then for all $k\in N$
we have
$$W^k(M)=W^k_0(M),\quad H^{2k}(M)=H^{2k}_0(M),\quad \mathrm{and}\quad
W^{2k}(M)=H^{2k}(M).$$
\end{lem}
\begin{proof} For the proof we refer to \cite{Sa}. The fact that 
$C^\infty_0(M)$ is dense in $H^{2k}(M)$ is an immediate consequence of 
\cite{Cn}. 
Indeed by \cite{Cn}, $(\Delta+\Id)^k$ is essentially self-adjoint on 
$C_0^\infty(M)$ for all $k\in\N$. Thus
\begin{equation}\label{3.3}
\overline{(\Delta+\Id)^k (C_c^\infty(M))}=L^2(M). 
\end{equation}
Let $f\in H^{2k}(M)$. Then $(\Delta+\Id)^kf\in L^2(M)$ and hence, by 
(\ref{3.3}) there exists a sequence $\{\varphi_j\}\subset C_c^\infty(M)$ such
that 
$$\parallel f-\varphi_j\parallel_{H^{2k}}=\parallel(\Delta+\Id)^{k}(f-\varphi_j)\parallel_{L^2}\to0$$
as $j\to\infty$.
\end{proof}

Under additional assumptions on $\xi$, similar results hold for weighted
Sobolev spaces \cite{Sa}. In general the following weaker results hold.

\begin{lem} \label{l3.3} For all $k\in\N$, the natural  inclusion 
$W_\xi^{2k}(M) \hookrightarrow H^{2k}_\xi(M)$ is bounded.
\end{lem}
\begin{proof} Let $k\in\N$.  Let $f\in W^{2k}_\xi(M)$. Then we have
$\nabla^jf\in L^2_\xi(M)$ for $j=0,...,2k$. Recall that 
\begin{equation}\label{3.5}
\Delta=-\Tr(\nabla^2f)
\end{equation}
and  $\nabla\Tr=0$. Hence it follows that there exists $C>0$ such that
$$|\Delta^jf|(x)\le C|\nabla^{2j}f|_g(x)$$
for all $j=0,...,k$ and $x\in M$. This implies $\Delta^jf\in L^2_\xi(M)$ for
$j=0,...,k$, and 
$$\parallel f\parallel_{H^{2k}_\xi}\le C\parallel f\parallel_{W^{2k}_\xi}.$$
\end{proof}

In order to deal with the inclusion in the other direction, we need some
preparation. Let $B_s\subset\R^n$ denote the ball of radius $s>0$ around
the origin in $\R^n$. Given $m\in\N$ and $r,K,\lambda>0$, denote by
$\Ell^m(r,K,\lambda)$ the set of elliptic differential operators
$$P=\sum_{|\alpha|\le m}a_\alpha(x)D^\alpha$$
of order $m$ in $B_r$ such
that the coefficients of $P$ satisfy:
\begin{enumerate}
\item $a_\alpha\in C^m(B_r)$.
\item$\sum_{|\alpha|<m}\p a_\alpha\p_{C^0(B_r)}\le K$, $\sum_{|\alpha|=m}\p
a_\alpha\p_{C^1(B_r)}\le K$.
\item$\lambda^{-1}\p \xi\p^m\le\sum_{|\alpha|=m}a_\alpha(x)\xi^\alpha\le
\lambda \p \xi\p^m$
for all $\xi\in\R^n$ and $x\in B_r$.
\end{enumerate}

Given an open subset $\Omega\subset \R^n$ and $k\in\N$,  $W^k(\Omega)$
is the usual Sobolev space.
\begin{lem}\label{l3.3a} 
Let $K,\lambda>0$ be given. There exists $r_0=r_0(K,\lambda)>0$
and $C=C(\lambda)>0$ such that for all $r\le r_0$, $P\in\Ell^m(r,K,\lambda)$ 
and $x_0\in B_r$:
$$\p u\p_{W^m(B_r)}\le C\left(\p Pu\p_{L^2(B_r)}+\p
  u\p_{L^2(B_r)}\right)$$
for all $u\in C_0^\infty(B_r)$
\end{lem}

\begin{proof}
Let $1\ge r>0$ and let $P\in\Ell^m(r,K,\lambda)$. Put 
$$P_0=\sum_{|\alpha|=m}a_\alpha(0)D^\alpha.$$
By Lemma 17.1.2 of \cite{H} there exists $C_1>0$ which depends only on
$\lambda$ such that for all $u\in C_0^\infty(B_r)$:
$$\p u\p_{W^m(B_r)}\le C\left(\p P_0u\p_{L^2(B_r)}+\p u\p_{L^2(B_r)}\right).$$
Now $Pu=P_0u+(P-P_0)u$. Thus
\begin{equation}\label{3.5a}
\p u\p_{W^m(B_r)}\le C\left(\p Pu\p_{L^2(B_r)}+\p (P-P_0)u\p_{L^2(B_r)}+\p
u\p_{L^2(B_r)}\right).
\end{equation} 
Next observe that
$$(P-P_0)u=\sum_{|\alpha|=m}(a_\alpha(x)-a_\alpha(0))D^\alpha
u+\sum_{|\alpha|<m}a_\alpha(x)D^\alpha u.$$
Hence by 2):
\begin{equation}\label{3.6}
\begin{split}
\p (P-P_0)u\p_{L^2(B_r)}\le& r\sum_{|\alpha|=m}\p a_\alpha\p_{C^1(B_r)}\p
u\p_{W^m(B_r)}\\
&\quad+\sum_{|\alpha|<m}\p a_\alpha\p_{C^0(B_r)}\p
u\p_{W^{m-1}(B_r)}\\
&\le K\left( r\p u\p_{W^m(B_r)}+\p u\p_{W^{m-1}(B_r)}\right).
\end{split}
\end{equation}
By the Poincar\'e inequality there exists $C_2>0$ which is independent of
$r\le 1$ such that for all $u\in C_0^\infty(B_r)$:
$$\p u\p_{W^{m-1}(B_r)}\le r\,C_2\p u\p_{W^m(B_r)}.$$
Using  this inequality, it follows from (\ref{3.6}) that
$$\p (P-P_0)u\p_{L^2(B_r)}\le  r\,C(K)\p u\p_{W^m(B_r)}.$$
Together with (\ref{3.5a}) we get
$$(1-rCC(K))\p u\p_{W^m(B_r)}\le C\left(\p Pu\p_{L^2(B_r)}+\p
  u\p_{L^2(B_r)}\right).$$ 
Set 
$$r_0=\min\{1,\frac{1}{2CC(K)}\}.$$
Then it follows that for all $r\le r_0$ and $u\in C^\infty_0(B_r)$:
$$\p u\p_{W^m(B_r)}\le 2C\left(\p Pu\p_{L^2(B_r)}+\p
  u\p_{L^2(B_r)}\right).$$ 
\end{proof}
\begin{lem} \label{lem3.5} Let $k\ge 1$ be even. Assume that $M$ has
 bounded curvature of order $k$. Let $K>0$ be such that 
$\sup_{x\in M}|\nabla^lR(x)|\le K$, $l=0,...,2k$.
There  exist constants $r_0=r_0(K)>0$ and $C=C(K)>0$ such that for all $x_0\in
M$ and $r\le \min\{r_0,\tilde\imath(x_0)\}$ one has
$$\parallel u\parallel_{W^{2k}(B_r(x_0))}\le C \parallel u
\parallel_{H^{2k}(B_r(x_0))}$$
for all $u\in C^\infty_0(B_r(x_0))$.
\end{lem}
\begin{proof} 
By \cite[Corollary 2.6 and 2.7]{Ei2} there exists a constant $C_1>0$,
which depends only on $K$, such that for every $x_0\in M$, every  
$r\le \tilde \imath(x_0)$, and  all $i,j,k=1,...,n$, one has
\begin{equation}\label{3.7}
\sup_{x\in B_r(x_0)}|D^\alpha g_{ij}(x)|\le C_1,\;\;|\alpha|\le 2k,\quad
\sup_{x\in B_r(x_0)}|D^\beta \Gamma^i_{jk}(x)|\le C_1,\;\; |\beta|\le 2k-1,
\end{equation}
where the $g_{ij}$ and $\Gamma^i_{jk}$ denote the 
coefficients of $g$ and
$\nabla$, respectively, with respect to normal coordinates on the geodesic ball
$B_r(x_0)$ of radius $r$ with center $x_0$. 

Let $x_0\in M$ and $r\le\tilde\imath(x_0)$. Let $B_r\subset T_{x_0}M$ 
denote the ball of radius $r$ around the origin. Let $W^{2k}(B_r)$ be the
Sobolev space with respect to the flat connection.  
Then it follows from (\ref{3.7}) that there exists $C_2=C_2(K)>0$ such that
\begin{equation}\label{3.8}
C_2^{-1}\parallel u\circ\exp_{x_0}\parallel_{W^{2k}(B_r)}\le 
\parallel u\parallel_{W^{2k}(B_r(x_0))}\le C_2\parallel u\circ\exp_{x_0}
\parallel_{W^{2k}(B_r)}
\end{equation}
for all $x_0\in M$, $r\le \tilde\imath(x_0)$, and $u\in C^\infty_c(B_r(x_0))$.
Let $\tilde g$ be the metric on $B_r$ which is the pull-back of
$g\upharpoonright B_r(x_0)$  with respect to
 $\exp_{x_0}\colon B_r\to B_r(x_0)$. Let $\widetilde\Delta$ be the 
Laplacian on $B_r$ with respect to $\tilde g$. Then by (\ref{3.8}) it is 
sufficient to show that there exists $C_3=C_3(K)>0$ such that
\begin{equation}\label{3.9a}
\parallel f\parallel_{W^{2k}(B_r)}\le C_3\parallel(\widetilde\Delta+\Id)^k f
\parallel_{L^2(B_r)}
\end{equation}
for all $x_0\in M$, $r\le \tilde\imath(x_0)$, and $f\in C^\infty_0(B_r)$. 
Set $P=(\widetilde\Delta+\Id)^{k}$. By
(\ref{3.7}) there exists $C_4>0$, which depends only on $K$, such that
$P\in\Ell^{2k}(r,1,C_4)$. Then by Lemma \ref{l3.3a},
there exist $r_0>0$ and $C_3>0$ such that (\ref{3.9a}) holds for all $x_0\in M$
and $r\le \min\{r_0,\tilde\imath(x_0)\}$. This completes the proof of the
lemma. 
\end{proof}

\begin{lem} \label{l3.6} Let $k\in\N$ be even. Suppose that $(M,g)$ has 
bounded curvature of order $2k$. Let $\beta: M\to\R^+$  be a function of
controlled decay. Then there
exists a canonical bounded inclusions 
\begin{equation} H^k_{\beta\tilde{\imath}^{-2kn}}(M) \hookrightarrow 
W^k_{\beta}(M)\quad\mathrm{and}\quad H^k_{\beta}(M) \hookrightarrow 
W^k_{\beta\tilde{\imath}^{2kn}}(M).
\end{equation}
\end{lem}
\begin{proof} By Theorem~\ref{coverthm}, there exists a covering 
$\{ B_{\frac{\tilde{\imath}}{2^k}(x_i)}(x_i) \}_{i=1}^\infty$ of $M$ by
balls and a constant $C>0$ such that 
\begin{equation}\label{3.9} 
\forall x\in M:\#\{x_i\mid x\in B_{\tilde{\imath}(x_i)}
(x_i) \} \leq C. 
\end{equation}

Let $\varphi\in C^\infty(\R)$ be such that $\varphi=1$ on 
$[0,1]$ and $\varphi=0$ on $[2,\infty)$. For $x \in M$ and 
$1\leq j \leq k$, we define
\begin{equation} \varphi_{j,x}(y) = \left\{ \begin{array}{ll} 
\varphi(2^j\frac{d(x,y)}{\tilde{\imath}(x)}), &  
y\in B_{\tilde{\imath}(x)}(x); \\ 0, & \hbox{otherwise}. \end{array}\right. 
\end{equation}
Then $\varphi_{j,x}\in C_0^\infty(M)$. Let $f\in H^k(M)$. Using Lemma 
\ref{l3.1}, it follows that $\varphi_{j,x}f\in H^k(B_{\tilde\imath(x)}(x))$. 
Then by
Lemma \ref{lem3.5} we get $\varphi_{j,x}f\in W^k(B_{\tilde\imath(x)}(x))
$ and
\[ \nabla^j(\varphi_{k,x}f) =  \sum_{p=0}^j\binom{j}{p}(\nabla^p
\varphi_{k,x})(\nabla^{j-p}f),\quad j=0,...,k. \]
By estimating the supremum-norm of the derivatives of $\varphi_{k,x}$ and 
using Lemma \ref{lem3.5} , we get
\begin{equation}
\begin{split}
 \|\varphi_{k,x}f\|_{W^k}& \le C\|f\|_{W^k \big(B_{\frac{\tilde{\imath}}
{2^{k-1}}(x)}(x)\big)}+C'\sum_{p=1}^k \binom{k}{p}\tilde{\imath}^{-p}(x)\|
\varphi_{k-1,x}f\|_{W^{k-p}} \\
&\leq C\|f\|_{H^k\big(B_{\frac{\tilde{\imath}}{2^{k-1}}(x)}(x)\big)}
+C''\sum_{p=1}^k \binom{k}{p}\tilde{\imath}^{-p}(x)
\|\varphi_{k-1,x}f\|_{H^{k-p}}. 
\end{split}
\end{equation}
By induction, this yields
\begin{equation}\label{3.15b} 
\|\varphi_{k,x_i}f\|_{W^k} \leq C\tilde{\imath}^{-k}(x_i)
\|f\|_{H^k(B_{\tilde{\imath}(x_i)}(x_i))}. 
\end{equation}
Let $f\in H^k_\beta$. By Lemma \ref{l1.4}, Lemma \ref{lem3.5}, (\ref{3.9})
and (\ref{3.15b}) we get
\begin{equation} \label{3.16}
\begin{split}
\|f\|_{W^k_{\tilde{\beta}}} \leq C\sum_{i=1}^\infty 
\beta^{\frac{1}{2}}(x_i)\|\varphi_{k,x_i}f\|_{W^k} \tilde{\imath}^{k}(x_i)
&\leq C\sum_{i=1}^\infty\beta^{\frac{1}{2}}(x_i) 
\|\varphi_{k,x_i}f\|_{H^k} \\ 
& \leq  C \sum_{i=1}^\infty \beta^{\frac{1}{2}}(x_i)
\tilde{\imath}^{-k}(x_i)\|f\|_{H^k({B_{\tilde{\imath}(x_i)}(x_i)})}. 
\end{split}
\end{equation}
By (\ref{2.2}) there exists a constant $C_1>0$ such that 
$$\tilde{\imath}(x_i)^{-k}\tilde{\imath}(x)^{kn}\le C_1$$ 
for all $i\in\N$ and $x\in B_{\tilde{\imath}(x_i)}(x_i)$. This implies
$$\sum_{i=1}^\infty \beta^{\frac{1}{2}}(x_i)
\tilde{\imath}^{-k}(x_i)\|f\|_{H^k({B_{\tilde{\imath}(x_i)}(x_i)})}\le C_2
\parallel f\parallel_{H^k_{\tilde{\imath}^{-2kn}\beta}},$$
which together with (\ref{3.16}) gives the first inclusion. The
proof of the second inclusion is analogous. 
\end{proof}

\begin{rem} Lemma \ref{l3.6} is not optimal. Under additional assumptions 
on $\beta$ one can show that $W^{2k}_\beta(M)=H^{2k}_\beta(M)$ \cite{Sa}.  
\end{rem} 

\section{Functions of the Laplacian.} \label{functions}
\setcounter{equation}{0}
\sectionmark{GEOMETRY OF MANIFOLDS WITH BOUNDED CURVATURE.}

Assume that $(M,g)$ is complete. Then $\Delta\colon C^\infty_c(M)\to L^2(M)$
is essentially self-adjoint and functions $f(\sqrt{\Delta})$ can be defined by
the spectral theorem for unbounded self-adjoint operators by
$$f(\sqrt{\Delta})=\int_0^\infty f(\lambda)dE_\lambda,$$
where $dE_\lambda$ is the projection spectral measure associated with 
$\sqrt{\Delta}$. Let $f\in L^1(\R)$ be even and let 
$$\hat f(\lambda)=\int_{-\infty}^\infty f(x)\cos(\lambda x)\;dx.$$
Then $f(\sqrt{\Delta})$ can also be defined by
\begin{equation}\label{4.1}
f(\sqrt{\Delta})=\frac{1}{2\pi}\int_{-\infty}^\infty \hat f(\lambda)
\cos(\lambda\sqrt{\Delta})\;d\lambda.
\end{equation}
This representation has been used in \cite{CGT}
to study the kernel of $f(\sqrt{\Delta})$. We will use (\ref{4.1}) to 
study $f(\sqrt{\Delta})$ as operator in weighted $L^2$-spaces. To this
end we  need  to study $\cos(\lambda\sqrt{\Delta})$ as operator in 
$L^2_\beta(M)$. Given $s>0$, let
 $\kappa(M,g,s)$ be the constant introduced in Definition \ref{kappadef}.

\begin{theo} \label{th4.1} Assume that $(M,g)$ has bounded curvature.
 Let $\beta$ be a function of moderate decay.  Then  $\cos({s\sD})$ extends
to a bounded operator in $L^2_\beta(M)$ for all $s\in \R$ and there exist
$C,c>0$ 
 such that  
\begin{equation}\label{3.14a}
\|\cos({s\sD}) \|_{L^2_\beta,L^2_\beta}\le C e^{c|s|}, \quad s\in\R.
\end{equation}
Moreover $\cos({s\sD})\colon L^2_\beta(M)\to L^2_\beta(M)$ is strongly
 continuous in $s$.
\end{theo}
\begin{proof} Let  $s>0$. Choose a sequence  $\{x_k\}_{k=1}^\infty\subset M$
which 
minimizes $\kappa(M,g;s)$. For $k\in\N$ let $P_k$ denote the multiplication 
by the characteristic function of 
$B_s(x_k)\setminus\bigcup_{i=0}^{k-1} B_s(x_j)$.
 Then each $P_k$ is an 
orthogonal projection in $L^2(M)$ and $L^2_\beta(M)$, respectively. Moreover
 the
projections satisfy
$P_kP_{k'}=0$
 for $k\not= k'$ and $\sum_{k=1}^\infty P_k =1$, where the series is strongly
 convergent. Obviously the
 image of $P_k$ consists of functions with support in $B_s(x_k)$.
Now recall that $\cos(t\sD)$ has unit propagation speed \cite[p.19]{CGT}, 
i.e., 
$$\supp \cos(s\sD)\delta_x\subset \overline{B_{|t|}(x)}$$
for all $x\in M$ and $t\in\R$. Let $f\in L^2(M)$. Then it follows that
$$\supp \cos(s\sD)P_kf\subset B_{2s}(x_k)$$
and 
$$\supp \cos(s\sD)\left((1-\chi_{B_{3s}(x_k)})f\right)\subset M- B_{2s}(x_k).$$
Hence
\begin{equation}\label{3.15}
\begin{split}
\parallel\cos(s\sD)f\parallel^2_\beta&=\sum_{k=1}^\infty\langle
\cos(s\sD)P_kf,\cos(s\sD)f\rangle_\beta\\
&=\sum_{k=1}^\infty\langle\cos(s\sD)P_kf,\cos(s\sD)
(\chi_{B_{3s}(x_k)}f)\rangle_\beta.
\end{split}
\end{equation}
Now observe that  the norm of $\cos(s\sD)$  as an operator in $L^2(M)$ is
bounded by $1$. This implies
\begin{equation*}
\begin{split}
\big|\langle\cos(s\sD)P_kf,\cos(s\sD)&(\chi_{B_{3s}(x_k)}f)\rangle_\beta\big|\\
&\le \sup_{y\in B_{3s}(x_k)}\beta(y)\parallel P_kf\parallel_{L^2}\cdot\parallel
\chi_{B_{3s}(x_k)}f\parallel_{L^2}.
\end{split}
\end{equation*}
To estimate the right hand side, we write 
\begin{equation*}
\sup_{y\in B_{3s}(x_k)}\beta(y)\parallel P_kf\parallel_{L^2}^2
=\int_M|P_kf(x)|^2\sup_{y\in B_{3s}(x_k)}
\left(\frac{\beta(y)}{\beta(x)}\right)\beta(x)\;dx.
\end{equation*}
Since the support of $P_kf$ is contained in $B_s(x_k)$, we can use (\ref{1.8})
to estimate the right hand side. This gives
$$\sup_{y\in B_{3s}(x_k)}\beta(y)\parallel P_kf\parallel_{L^2}^2
\le C_\beta^{-1}\frac{1}{\beta(1+4s)}\parallel P_kf\parallel_{L^2_\beta}^2.$$
A similar inequality holds with respect to
$\parallel\chi_{B_{3s}(x_k)}f\parallel_{L^2}$. Putting the estimations 
together, we get
\begin{equation}\label{3.15a}
\begin{split}
\big|\langle\cos(s\sD)P_kf,\cos(s\sD)&(\chi_{B_{3s}(x_k)}f)\rangle_\beta\big|\\
&\le C_\beta^{-1}\frac{1}{\beta(1+6s)}\parallel P_kf\parallel_{L^2_\beta}
\cdot \parallel\chi_{B_{3s}(x_k)}f\parallel_{L^2_\beta}.
\end{split}
\end{equation}
Now recall that by Lemma \ref{lordjagged} we have  $\kappa(M,g;s)<\infty$.
Hence we get
\begin{equation*} 
\sum_{k=1}^\infty \|\chi_{B_{3s}(x_k)}f\|^2_{L^2_\beta} \leq 
\kappa(M,g;s)\|f\|^2_{L^2_\beta} < \infty. 
\end{equation*}
Together with \ref{3.15} and (\ref{3.15a}) we obtain
\begin{equation}\label{3.17}
\begin{split}
\parallel\cos(s\sD)f\parallel^2_{L^2_\beta}&\le
C_\beta^{-1}\frac{1}{\beta(1+6s)} \parallel f\parallel_{L^2_\beta}
\sum_{k=1}^\infty\parallel \chi_{B_{3s}(x_k)}f\parallel_{L^2_\beta}\\
&\le C_\beta^{-1}\frac{1}{\beta(1+6s)}\kappa(M,g,s)^{1/2}\parallel
f\parallel_{L^2_\beta}^2. 
\end{split}
\end{equation}
Recall that by (1.10) we have  $\beta(x)\le C(1+d(x,p))^{-1}$, $x\in
M$. Therefore $L^2(M)\subset L^2_\beta(M)$, and 
 $L^2(M)$ is a dense subspace of $L^2_\beta(M)$. This implies
 that $\cos(s\sD)$ extends to a bounded operator in
$L^2_\beta(M)$. Moreover by (\ref{1.7a}) and Lemma \ref{lordjagged} it 
follows that there exist constants $C,c>0$ such that
$$\parallel \cos(s\sD)\parallel^2_{L^2_\beta,L^2_\beta}\le C e^{cs},\quad  
s\in [0,\infty).$$ 
Since $\cos(-s\sD)=\cos(s\sD)$, this extends to all $s\in\R$ such that 
(\ref{3.14a}) holds. The strong continuity is a consequence of the local
bound of the norm and the strong continuity on the dense subspace 
$L^2 (M)\subseteq L^2_\beta(M)$.
\end{proof}

Using  Theorem~\ref{th4.1}, we can study $f(\sqrt{\Delta})$
as an operator in $L^2_\beta(M)$. Given $c\ge0$, let
$$\mathcal{F}^1(c)=\bigg\{f\in L^1(\R)\colon \int_{-\infty}^\infty |\hat{f}
(\lambda)|e^{c|\lambda|} d\lambda < \infty\bigg\}.$$

\begin{lem} \label{l4.2} Assume $(M,g)$ has bounded curvature and
let $\beta $ be a function of moderate decay. Then there exists a
 constant $c= c(M,g,\beta)$, such that for all even functions
 $f\in\mathcal{F}^1(c)$, the operator $f(\sD)$ extends to
 a bounded operator 
in $L^2_\beta(M)$. 
Moreover, there exists a constant $C_1=C_1(M,g,\beta)>0$ such that
\begin{equation} \|f(\sD)\|_{L^2_\beta,L^2_\beta} \leq C_1
\|\hat{f}\|_{L^1_{e^{c|\cdot|}}} \end{equation}
for all $f$ as above.
If $\kappa(M,g;s)$ is at most sub-exponentially increasing, then 
$c(M,g,\beta)>0$ can be chosen arbitrarily.
\end{lem}
\begin{proof} By Theorem \ref{th4.1} there exist constants
$C,c>0$, depending on $(M,g,\beta)$,  such that 
$$ \|\cos(\sD)\|_{L^2_\beta,L^2_\beta}\le C e^{c|s|},$$ 
for all  $s\in\R$.  Let $\varphi\in L^2(M)$. Using (\ref{4.1}), it follows
that 
\begin{equation}\label{4.7} 
\|f(\sD)\varphi\|_{L^2_\beta} \le\frac{C}{\sqrt{2\pi}}
\|\hat{f}\|_{L^1_{e^{C|\cdot|}}}\|\varphi\|_{L^2_\beta}. 
\end{equation}
Since $L^2(M)$ is dense in $L^2_\beta(M)$, it follows from (\ref{4.7}) 
that $f(\sD)$ extends to a bounded operator in
$L^2_\beta(M)$. The last statement is obvious. \end{proof}

\begin{rem} It is not difficult to see, that (\ref{4.1}) is in
  fact strongly convergent in $L^2_\beta$. \end{rem} 
 
\begin{cor} \label{c4.3} Suppose that $(M,g)$ has bounded
  curvature and let $\beta$ be a function of moderate decay. Then the
  following holds:

\noindent
a) For every $t>0$, the heat  operator $e^{-t\Delta}$ extends to bounded
  operator in $L^2_\beta(M)$. Its norm is uniformly bounded in $t$ on compact
  intervals of $\R^+$.

\noindent  
b) In the region $\left\{\lambda\in\C\colon\Re(\sqrt{-\lambda}) > c(M,g,\beta)
\right\}$ the resolvent
  $(\Delta-\lambda)^{-1}$ extends to a bounded operator in $L^2_\beta(M)$.
  The  function of $\lambda\mapsto (\Delta-\lambda)^{-1}$ is locally bounded
 and holomorphic on this domain.

\noindent
c) If $\beta$ is of sub-exponential decay and $\kappa(M,g;s)$
  is  at most sub-exponentially increasing for $s>s_0$, then
  $(\Delta-\lambda)^{-1}:L^2_\beta(M)\mapsto L^2_\beta(M)$
 is defined and bounded for all $\lambda \in \C\setminus [0,\infty)$. 

\end{cor}
\begin{proof} This follows from Lemma~\ref{l4.2} and

\begin{equation} \int_{-\infty}^\infty 
e^{-t x^2}\cos(xy)dx = \sqrt{\frac{\pi}{t}}e^{-\frac{y^2}{4t}}, 
\quad \int_{-\infty}^\infty 
\frac{1}{\lambda+x^2}\cos(xy)dx
=\frac{\pi}{\sqrt{\lambda}}e^{-\sqrt{\lambda}\,|y|}.\end{equation}
\end{proof}

Let $\beta$ be of moderate decay. 
There is a canonical pairing $(\cdot,\cdot)$ between 
$L^2_\beta(M)$ and $L^2_{\beta^{-1}}(M)$ given by
$$(f,g)=\int_M f(x)g(x)\;dx, \quad f\in L^2_\beta(M),\; 
g\in L^2_{\beta^{-1}}(M).$$
This pairing is non-degenerate so that $L^2_{\beta^{-1}}(M)$ is canonically
isomorphic to the dual of $L^2_\beta(M)$. Moreover, we have the following
inclusions
$$L^2_{\beta^{-1}}(M)\subset L^2(M)\subset L^2_\beta(M).$$
By duality it follows that Theorem \ref{th4.1}, Lemma \ref{l4.2} and Corollary 
\ref{c4.3} also hold w.r.t. $\beta^{-1}$. 
Especially, it follows  that $f(\sD)$ defined on $L^2_{\beta^{-1}}(M)$ is the
 restriction of $f(\sD)_{|L^2}$. Moreover, we have the identity
\begin{equation} f(\sD)_{|L^2_{\beta^{-1}}} = 
\left(\bar{f}(\sD)_{|L^2_\beta}\right)^\ast. \end{equation}



\begin{lem} \label{l4.5} Let $\beta$ be a function of moderate  
decay. If $\lambda$ and $\bar \lambda$ satisfy condition b) of Corollary
\ref{c4.3}, then
$$H^2_\beta(M)=(\Delta-\lambda)^{-1}(L^2_\beta(M)).$$
\end{lem}
\begin{proof} First note that $C_0^\infty(M)$ is dense in
  $L^2_\beta(M)$. Indeed  $C_0^\infty(M)$ is dense in $L^2(M)$, and 
$L^2(M)$ is dense in $L^2_\beta(M)$.
Let $f=(\Delta-\lambda)^{-1}g$, $g\in L^2_\beta(M)$. Then there exists a
sequence 
$\{\varphi_i\}_{i\in\N}\subset C^\infty_0(M)$ which converges to $g$ in
$L^2_\beta(M)$  and
$(\Delta-\lambda)^{-1}\varphi_i$ converges to $f$ in $L^2(M)$. Let $\varphi\in
C^\infty_0(M)$. Then
$$\langle f,\Delta\varphi\rangle=\lim_{i\to\infty}\langle(\Delta-\lambda)^{-1}
\varphi_i,\Delta\varphi\rangle=\lim_{i\to\infty}\langle \varphi_i+\lambda
(\Delta-\lambda)^{-1}\varphi_i,\varphi\rangle=\langle g+\lambda f,\varphi
\rangle.$$
Thus $\Delta f=g+\lambda f\in L^2_\beta(M)$ and hence $f\in
H^2_\beta(M)$. Now suppose that 
$f\in H^2_\beta(M)$ and set $g=(\Delta-\lambda)f$. Then $g\in L^2_\beta(M)$
and  
we need to show that $f=(\Delta-\lambda)^{-1}g$. Let $\varphi\in
C_0^\infty(M)$. By definition of $(\Delta-\lambda)^{-1}g$, there exists a
sequence $\{g_i\}_{i\in\N}\subset L^2(M)$ such that $(\Delta-\lambda)^{-1}g_i$ 
converges to $(\Delta-\lambda)^{-1}g$ in $L^2_\beta(M)$ as $i\to\infty$. Using
this fact, we get
\begin{equation}\label{4.10}
\langle (\Delta-\lambda)^{-1}g,\varphi\rangle=\langle g,
(\Delta-\bar\lambda)^{-1}\varphi\rangle =\langle(\Delta-\lambda)f,
(\Delta-\bar \lambda)^{-1}\varphi\rangle.
\end{equation}
Now observe that $(\Delta-\bar \lambda)^{-1}\varphi$ belongs to $H^2(M)$. By
Lemma \ref{l3.1}, there exists a sequence $\{\varphi_i\}_{i\in\N}\subset 
C_0^\infty(M)$ which converges to $(\Delta-\bar \lambda)^{-1}\varphi$ in
$H^2(M)$. Thus
$$\langle(\Delta-\lambda)f,(\Delta-\bar \lambda)^{-1}\varphi\rangle
=\lim_{i\to\infty}\langle (\Delta-\lambda)f,\varphi_i\rangle=
\langle f,(\Delta-\bar\lambda)\varphi_i\rangle=\langle f, \varphi\rangle.$$
Together with (\ref{4.10}) this implies that $f=(\Delta-\lambda)^{-1}g$.
\end{proof}

\begin{lem}\label{l4.6}
Let $\beta$ be a function of moderate decay. Then
$(\Delta+\lambda)(C^\infty_0(M))$ is 
dense in $L^2_\beta(M)$ for every $\lambda\in \R^+$. 
\end{lem}
\begin{proof}
As in (\ref{3.3}) it follows from the essential self-adjointness of 
$\Delta+\lambda\Id$ that $(\Delta+\lambda)(C_c^\infty(M))$ is dense in
$L^2(M)$.  
Moreover since $\beta$ is monotonically decreasing, we have that $L^2(M)\subset
L^2_\beta(M)$ is dense and $\parallel f\parallel_\beta\le C\parallel f\parallel$
for $f\in L^2(M)$.
This implies that $(\Delta+\lambda)(C_c^\infty(M))$ is also dense in 
$L^2_\beta(M)$.
\end{proof}
\begin{cor} \label{c4.7}
Let $\beta$ be of moderate decay. 
Then $C_0^\infty(M)$ is dense in $H^2_\beta(M)$.
\end{cor}
\begin{proof} 
Let $f\in H^2_\beta(M)$. Let $\lambda\gg 0$.
By  Lemma \ref{l4.5} there exists $g\in L^2_\beta(M)$ 
such that $f=(\Delta+\lambda)^{-1}g$. By Lemma \ref{l4.6} there exists a 
sequence $\{\varphi_i\}_{i\in\N}\subset C_c^\infty(M)$ such that 
$(\Delta+\lambda)\varphi_i$ converges to $g$ in $L^2_\beta(M)$ as $i\to\infty$.
Thus $\varphi_i\to f$ in $L^2_\beta(M)$ and $(\Delta+\lambda)\varphi_i$ 
converges to $g=(\Delta+\lambda)f$ as $i\to\infty$. This implies that 
$\varphi_i$ converges to $f$ in $H^2_\beta(M)$.
\end{proof}

\section{Equivalent Metrics and Sobolev Spaces.}
\sectionmark{Equivalent Metrics and Sobolev Spaces.}
\setcounter{equation}{0}

In this section we study the dependence of the Sobolev spaces on the
metric. We will  prove, that if $g\sim_\beta^k h$ for 
an appropriate $\beta$, then the  Sobolev spaces defined with 
respect to 
$g$ and $h$ are equivalent up to order $k$.
We assume that all metrics  have bounded sectional curvature. To indicate the
dependence of the corresponding 
Sobolev space on the Riemannian metric $g$, we will write $W^k_\xi(M;g)$ and 
$H^{2k}_\xi(M;g)$, respectively.

\begin{lem}\label{l5.1}  Let $\beta$ be of moderate
decay. Assume that $g\sim_\beta^k h$. Then the Sobolev 
spaces $W^{k}_\xi(M;g)$ and $W^{k}_\xi(M;h)$ are equivalent.
\end{lem}
\begin{proof} First note that by Lemma \ref{l1.7} the metrics $g$ and $h$ 
are quasi-isometric. This implies that $L^2_\xi(M,g)$ and 
$L^2_\xi(M;h)$ are equivalent. So the statement of the lemma holds for $k=0$. 
Let $f\in C^\infty(M)$. Let $k\ge 1$.  By induction we will prove that for
$l\le k$ there
 exists $C_l>0$ such that for $a,b\in\N_0$, $a+b=l$,  
\begin{equation}\label{5.1}
\big|(\nabla^g)^a(\nabla^h)^b f\big|_h(x) \leq 
C_l\sum_{i=0}^{a+b} \big|(\nabla^g)^if|_g(x), \quad x\in M.
\end{equation}
Let $l=1$. Since on functions the connections equal $d$, 
(\ref{5.1}) follows from quasi-isometry of  $g$ and $h$.

Next suppose that (\ref{5.1}) holds for $1\le l<k$. To establish (\ref{5.1})
for  
$l+1$, we proceed by induction with respect to $a$. Let $a,b\in\N$ with
$a+b=l+1$. If $a=l+1$ there is nothing to prove. Let $a<l+1$. Then
\begin{equation} (\nabla^g)^a(\nabla^h)^bf = (\nabla^g)^a(\nabla^h-\nabla^g)
(\nabla^h)^{b-1}f+(\nabla^g)^{a+1}(\nabla^h)^{b-1}f.
\end{equation}
and therefore, we get
\begin{equation*} 
\begin{split}
|(\nabla^g)^a(\nabla^h)^bf|_h(x) \le 
|(\nabla^g)^a(\nabla^h&-\nabla^g)(\nabla^h)^{b-1}f|_h(x)\\
&+|(\nabla^g)^{a+1}(\nabla^h)^{b-1}f|_h(x),\quad x\in M.
\end{split}
\end{equation*}
Using $g\sim^k_\beta h$ together with the binomial formula and the
induction hypothesis, it follow
that (\ref{5.1}) holds for $l+1$. Especially, putting $a=0$ we get
\begin{equation}\label{5.3}
|(\nabla^h)^lf|_h(x) \leq C_l\sum_{i=0}^l |(\nabla^g)^if|_g(x),\quad x\in M,\;
l\le k.
\end{equation}
Suppose that $f\in C^\infty(M)\cap W^k_\xi(M;g)$. Then (\ref{5.3}) implies that
$f\in C^\infty(M)\cap W^k_\xi(M;h)$ and
$$\parallel f\parallel_{W^k_\xi(M;h)}\le C 
\parallel f\parallel_{W^k_\xi(M;g)}.$$
By Lemma \ref{l3.1}, $C^\infty(M)\cap W^k_\xi(M;g)$ is dense in
 $W^k_\xi(M;g)$. Therefore this inequality holds for all $f\in W^k_\xi(M,g)$.
By symmetry, a similar inequality holds with the roles of $g$ and $h$ 
interchanged. This concludes the proof.
\end{proof}

Next we compare the Sobolev spaces $H^{2k}_\xi(M;g)$ and $H^{2k}_\xi(M;h)$.
Let $\Delta_g$ denote the Laplace operator with respect to the metric $g$.
 Recall, that
\begin{equation*}
\Delta_g = (\nabla^g)^\ast\nabla^g, 
\end{equation*}
and that the formal adjoint $(\nabla^g)^\ast$ of $\nabla^g$ is given by
\begin{equation}\label{5.3a}
(\nabla^g)^\ast = -\Tr (g^{-1}\nabla^g),
\end{equation}
where $g^{-1}:T^*M\to TM$ is the isomorphism induced by the metric and
$\Tr \colon T^*M\otimes TM\to \R$ denotes the contraction.
Since $\nabla^g \Tr = 0$ and $\nabla^g g^{-1}=0$, we get
\begin{equation}\label{5.4}
 \Delta_g^k = (-1)^k(\Tr g^{-1})^k (\nabla^g)^{2k}.
\end{equation}

\begin{lem} \label{l5.2}
 Assume that $g\sim^{2k}_\beta h$. 
Then for each $l$, $0\leq l\leq 2k$ and $j$, $0\le j\le 2l$, there exist
sections $\xi_{jl}^g,\xi_{jl}^h\in C^\infty(\Hom((T^*M)^{\otimes j},\R))$ such
that 
\begin{equation}\label{5.4a}
\Delta_g^l-\Delta_h^l=\sum_{j=0}^{2l}\xi_{jl}^g\circ(\nabla^g)^j
=\sum_{j=0}^{2l}\xi_{jl}^h\circ(\nabla^h)^j
\end{equation}
and there exists $C>0$ such that for $0\le p\le l$ 
\begin{equation}\label{5.4b}
|(\nabla^g)^p\xi_{jl}^g|_g(x)\le C\beta(x),\quad 
|(\nabla^h)^p\xi_{jl}^h|_h(x)\le C\beta(x),\;x\in M.
\end{equation}
\end{lem}
\begin{proof} 

Using (\ref{5.4}) we get
\begin{equation}\label{5.5}
\begin{split}
(-1)^l(\Delta^l_g-\Delta^l_h) &= (\Tr g^{-1})^l(\nabla^g)^{2l} - 
(\Tr h^{-1})^l(\nabla^h)^{2l} \\
&= (\Tr g^{-1})^l\left((\nabla^g)^{2l}-(\nabla^h)^{2l}\right)+ 
\left( (\Tr g^{-1})^l-(\Tr h^{-1})^l\right)(\nabla^h)^{2l}. 
\end{split}
\end{equation}
First consider the second term. Note that there exists $C>0$ such that
\begin{equation}\label{5.6}
\begin{split}
\big|(\nabla^g)^p\left( (\Tr g^{-1})^l-(\Tr h^{-1})^l\right)|_g(x)
\le C |(\nabla^g)^p(g-h)|_g(x).
\end{split}
\end{equation}
Since $g\sim^{2k}_\beta h$, the right hand side is bounded by $C_1\beta(x)$. 
By symmetry, the same estimation holds with respect to $h$.

To deal with the first term on the right hand side of (\ref{5.5}), we use
$$(\nabla^g)^j-(\nabla^h)^j=(\nabla^g)^{j-1}(\nabla^g-\nabla^h)+
\left((\nabla^g)^{j-1}-(\nabla^h)^{j-1}\right)\nabla^h$$
and proceed by induction with respect to $j$.
\end{proof}

\begin{cor} \label{c5.3} 
Let $\beta$ be of controlled decay. Assume that 
$\beta \tilde{i}^{-2kn}$ is bounded, $g\sim^{2k}_\beta h$ and $(M,g)$ and $(M,h)$ have both bounded curvature of order $2k$. Then $H^{2k}_\rho(M,g)$ and
$H^{2k}_\rho(M,h)$ are equivalent for all functions $\rho$ of controlled decay.
\end{cor} 
\begin{proof}  
Let $f\in C^\infty(M)\cap H^{2k}_\rho(M;g)$. Using Lemma \ref{l3.6} and Lemma
\ref{l5.1} we get
$$\parallel f\parallel_{H^{2k}_\rho(M;g)}\ge 
C_1 \parallel f\parallel_{W^{2k}_{\tilde{\imath}^{4k}\rho}(M;g)}
\ge C_2 \parallel f\parallel_{W^{2k}_{\tilde{\imath}^{4kn}\rho}(M;h)}
\ge C_3 \parallel f\parallel_{W^{2k}_{\beta^2\rho}(M;h)}.$$
By Lemma \ref{l5.2} it follows that $f\in C^\infty(M)\cap H^{2k}_\rho(M,h)$
and there exists a constant $C>0$, which is independent of $f$,  such that
$$\parallel f\parallel_{H^{2k}_\rho(M;h)}\le C
\parallel f\parallel_{H^{2k}_\rho(M;g)}.$$ 
By symmetry, a similar inequality holds with $g$ and $h$ interchanged. 
\end{proof}

\section{Trace class estimates}  
\setcounter{equation}{0}

Let $(M,g)$ be an $n$-dimensional  Riemannian manifold with bounded sectional
curvature, $|K_M|<K$.  Let $e^{-t\Delta_g}(x,y)$ denote the heat kernel of
$\Delta_g$.  Let $0<a_1<a_2<\infty$. Let $\tilde{\imath}$ be the modified
injectivity radius defined by  (\ref{2.1}).
It follows from  \cite[Proposition
1.3]{CGT}, that there exist $C_1,c_1>0$ such that  
\begin{equation} \label{6.1}
e^{-t\Delta_g}(x,y) \le
C_1\tilde{\imath}(x)^{-\frac{n}{2}}
\tilde{\imath}(y)^{-\frac{n}{2}}e^{-c_1d^2(x,y)},\quad t\in[a_1,a_2].  
\end{equation}
Let $c<c_1$. Then by (\ref{6.1}) and (\ref{2.2}) there exists $C>0$ such that
\begin{equation} \label{6.2}
e^{-t\Delta_g}(x,y) \leq
C\tilde{\imath}(x)^{-\frac{n(n+1)}{2}}e^{-cd^2(x,y)},\quad  t\in[a_1,a_2].  
\end{equation}

\begin{lem} \label{lem6.1} Let $\beta$ be a function of moderate
  decay. Assume that there exist real numbers $a,b$ such that 
\begin{itemize} 
\item[i)] $a+b=2$,
\item[ii)] $\beta^b \in L^1(M)$,
\item[iii)] $\beta^a\tilde{\imath}^{-\frac{n(n+1)}{2}} \in L^\infty(M)$.
\end{itemize}
Let $M_\beta$ the operator of multiplication by $\beta$. Then for every  
$p\in\N_0$, the operator $M_\beta \Delta_g^pe^{-t\Delta_g}$ is
Hilbert-Schmidt. For $t$ in a compact interval in $\R^+$, the Hilbert-Schmidt
norm is bounded. 
\end{lem}
\begin{proof} 
we have
\begin{equation}\label{6.4a} 
M_\beta \Delta^pe^{-t\Delta} = \left(M_\beta
 e^{-\frac{t}{2}\Delta}\right)\left(\Delta^pe^{-\frac{t}{2}\Delta}\right). 
\end{equation}
Note that the operator norm of $\Delta^p e^{-\frac{t}{2}\Delta}$ is 
bounded on compact subsets of $\R^+$. Hence we may assume that $p=0$. By
Corollary 
\ref{c4.3}, 1), it follows that $e^{-t\Delta}$ extends to a bounded operator
in $L^2_{\beta^b} (M)$ and its norm is uniformly bounded for $0<a\le t\le b$. 
The condition $\beta^b \in L^1(M)$ implies that $1 \in L^2_{\beta^b}$. Hence
$e^{-t\Delta}1\in L^2_{\beta^b}(M)$. Let $e^{-t\Delta}(x,y)$ be the kernel of
$e^{-t\Delta}$. Then
$$\langle 1,e^{-t\Delta}1\rangle_{L^2_{\beta^b}}=\int_M \int_M
\beta^b(x)e^{-t\Delta_g}(x,y) dydx.$$  
The integral converges since $e^{-t\Delta}(x,y)$ is positive.
Thus we get
\begin{equation*}
\begin{split} 
\int_M\int_M |\beta(x)&e^{-t\Delta_g}(x,y)|^2dydx = \int_M\int_M
\beta^2(x)e^{-2t\Delta_g}(x,y) dydx\\  
& \leq \sup_{z,w\in M}|\beta^a(z)e^{-t\Delta_g}(z,w)| \int_M\int_M
\beta^b(x)e^{-t\Delta_g}(x,y) dydx \\ 
 &\leq C \sup_{z\in M}|\beta^a(z)\tilde{\imath}^{-\frac{n(n+1)}{2}}(z)|
 \int_M\beta^b(x)\left(e^{-t\Delta}(1)\right)(x) dx \\  
& \le  C_1 \parallel e^{-t\Delta}(1)\parallel_{L^2_{\beta^b}}. 
\end{split}
\end{equation*}
This proves the lemma. 
\end{proof}

\begin{lem}\label{lem6.2} Assume $\beta$ is a function of moderate decay
  and that there exist real numbers $a,b$ such that 
\begin{itemize} 
\item[i)] $b\geq 1$ and $a+b=2$,
\item[ii)] $\beta^\frac{b}{3} \in L^1(M)$,
\item[iii)] $\beta^{\frac{a}{3}}\tilde{\imath}^{-\frac{n(n+2)}{2}} \in
  L^\infty(M)$. 
\end{itemize}
Let $M_\beta$ be the operator of multiplication by $\beta$. Then the operator
$M_{\tilde{\imath}^{-2n}}M_\beta \Delta^pe^{-t\Delta}$ is a trace-class
operator for $p\in \N$. For $t$ in a compact interval, the trace-class norm is
bounded. 
\end{lem} 
\begin{proof} We decompose the operator as 
\begin{equation} \label{6.4} 
M_{\tilde{\imath}^{-2n}}M_\beta \Delta^pe^{-t\Delta} =
\left\{ M_{\tilde{\imath}^{-2n}}M_\beta e^{-\frac{t}{2}\Delta}
M_{\beta^{-\frac{1}{3}}}\right\}\cdot\left\{M_{\beta^{\frac{1}{3}}}
\Delta^pe^{-\frac{t}{2}\Delta^2}\right\}. 
\end{equation} 
Since $\beta$ is non-increasing and $\beta(x)\le 1/2$ outside a compact set,
it follows that $\beta^{\frac{1}{3}}\le C\beta^{\frac{b}{3}}$ for $b\ge 1$. 
Hence by ii) we get $\beta^{\frac{1}{3}}\in L^1(M)$. Moreover by iii) it
follows that $\beta^{\frac{a}{3}}\tilde{\imath}^{-\frac{n(n+1)}{2}} \in
 L^\infty(M)$. Hence by Lemma \ref{lem6.1}, the second factor on the right
hand side of (\ref{6.4}) is a Hilbert-Schmidt operator and its Hilbert-Schmidt
norm is bounded for $t$ in a compact interval in $\R^+$. It remains to show
that the first factor is Hilbert-Schmidt and that the Hilbert-Schmidt norm is
bounded on compact intervals.  By iii) we have
$$\beta^a\tilde{\imath}^{-\frac{n(n+1)}{2}-2n}\in L^\infty(M).$$
Using this observation together with (\ref{6.2}), we get
\begin{equation} 
\begin{split}
\int_M\int_M
  |\tilde{\imath}^{-2n}(x)&\beta(x)e^{-t\Delta}(x,y)\beta^{-\frac{1}{3}}(y)|^2
  dx dy \\
&\leq C\sup_{z\in
  M}|\tilde{\imath}^{-\frac{n(n+1)}{2}-2n}(z)\beta^a(z)| \int_M\int_M
  \beta^b(x)e^{-t\Delta}(x,y)\beta^{-\frac{2}{3}}(y) dx dy. 
\end{split}
\end{equation} 
Now observe that by ii), $\beta^{-\frac{2}{3}}$ belongs to
$L^2_{\beta^{\frac{b+4}{3}}}(M)$. Since $\beta^{\frac{b+4}{3}}\le
  C\beta^{\frac{b}{3}}$, it follows from ii) that $\beta^{\frac{b+4}{3}}$
is integrable. Hence by Corollary \ref{c4.3}, $e^{-t\Delta}$ extends to
a bounded operator in $L^2_{\beta^{\frac{b+4}{3}}}(M)$. Therefore
$\int_M e^{-\Delta}(x,y)\beta^{-\frac{2}{3}}(y) dy \in
L^2_{\beta^\frac{b+4}{3}},$ 
and the norm is uniformly bounded for $t$ in a compact interval of $\R^+$.
 Next note that
$\beta^b \in L^2_{\beta^{-\frac{b+4}{3}}}$. Hence  
\begin{equation} 
\int_M\int_M
  \beta^b(x)e^{-t\Delta}(x,y)\beta^{-\frac{2}{3}}(y) dx dy=\left<
  \beta^b,e^{-t\Delta}\beta^{-\frac{2}{3}} 
  \right><\infty. 
\end{equation}
This implies the lemma. 
\end{proof} 

\begin{lem} \label{lem6.3} 
Let $\beta$ be a function of moderate decay, satisfying the conditions of
Lemma \ref{lem6.2}. Let $g,h$ be two complete metrics on $M$ such that
$g\sim^2_\beta h$.  Let $\Delta_g$ and
$\Delta_h$ be the Laplacians of $g$ and $h$, respectively. Then 
$$(\Delta_g-\Delta_h)e^{-t\Delta_g}\quad\mathrm{and}\quad
e^{-t\Delta_g}(\Delta_g-\Delta_h)$$
are trace class operators, and the trace norm is uniformly bounded for $t$ in a
compact subset of $(0,\infty)$.  
\end{lem}
\begin{proof} 
We decompose $e^{-t\Delta_g}$ as
\begin{equation} 
e^{-t\Delta_g} = \left(
e^{-\frac{t}{2}\Delta_g}M_{\beta^{-\frac{1}{3}}}\right)\cdot\left(
M_{\beta^{\frac{1}{3}}}e^{-\frac{t}{2}\Delta_g}\right). 
\end{equation}
By Lemma \ref{lem6.1}, the second factor is a Hilbert-Schmidt operator and it
suffices to show that
$ (\Delta_g-\Delta_h)e^{-t\Delta_g}M_{\beta^{-\frac{1}{3}}}$
is Hilbert-Schmidt and 
that the Hilbert-Schmidt norm is bounded for $t$ in a compact interval. Using
Lemma \ref{l5.2} and Lemma \ref{l3.6},it follows that the Hilbert-Schmidt norm
can be estimated by
\begin{equation*}
\begin{split}
\parallel
(\Delta_g-\Delta_h)e^{-t\Delta_g}M_{\beta^{-\frac{1}{3}}}\parallel^2_2&\le
C\sum_{i=0}^2\int_M\int_M |(\nabla^g)^i e^{-t\Delta_g}(x,y)
\beta^{-\frac{1}{3}}(y)|^2_g \beta^2(x)\;dx dy \\
&=C\int_M\parallel e^{-t\Delta_g}(\cdot,y)
\beta^{-\frac{1}{3}}(y)\parallel^2_{W^2_{\beta^2}}\;dy\\
&\le C_1\int_M\parallel e^{-t\Delta_g}(\cdot,y)
\beta^{-\frac{1}{3}}(y)\parallel^2_{H^2_{\beta^2\tilde\imath^{-4n}}}\;dy\\
&\leq C_2\sum_{q=0}^1 \int_M \|\beta(\cdot)\tilde{\imath}^{-2n}(\cdot)
\Delta_g^q e^{-t\Delta_g}(\cdot,y)\beta^{-\frac{1}{3}}(y)\|^2_2 dy \\ 
&= C_2\sum_{q=0}^1 \|M_{\beta}M_{\tilde{\imath}^{-2n}} \Delta_g^q
e^{-t\Delta_g}M_{\beta^{-\frac{1}{3}}}\|^2_2. 
\end{split}
\end{equation*}
By Lemma \ref{lem6.2} the right hand side  is finite and bounded for $t$ in a
compact interval of $\R^+$. To prove that $e^{-t\Delta_g}(\Delta_g-\Delta_h)$
is a trace class operator, it suffices to establish it for its adjoint
 $(\Delta_g-(\Delta_h)^{*_g})e^{-t\Delta_g}$ with respect to $g$. By
 (\ref{5.4a}) and  (\ref{5.3a}) we have
\begin{equation}
\Delta_g-(\Delta_h)^{*_g}
=(\xi_{01}^g)^{*_g}+(\nabla^g)^{*_g}\circ(\xi_{11}^g)^{*_g}
+\big[(\nabla^g)^{*_g}\big]^2\circ(\xi_{21}^g)^{*_g}. 
\end{equation}
Using (\ref{5.3a}) and (\ref{5.4b}), it  follows that there exist
$\eta_j\in C^\infty(\Hom((T^*M)^{\otimes j},\R))$ such that
\begin{equation}\label{6.9a}
\Delta_g-(\Delta_h)^{*_g}=\eta_0+\eta_1\circ
\nabla^g+\eta_2\circ(\nabla^g)^2
\end{equation}
and these section satisfy
\begin{equation}\label{6.10}
|\eta_j|_g(x)\le C\beta(x),\quad 0\le j\le 2, \; x\in M.
\end{equation}
Using (\ref{6.9a}) and (\ref{6.10}) we can proceed as above and prove that
$(\Delta_g-(\Delta_h)^{*_g})e^{-t\Delta_g}$ is a trace class operator.
\end{proof}
We are now ready to prove Theorem \ref{th0.1}. We note that for
equivalent metrics, the Hilbert spaces $L^2(M,g)$ and $L^2(M,h)$ are
equivalent. Hence we may regard $e^{-t\Delta_h}$ as bounded operator in
$L^2(M,g)$.

{\bf Proof of Theorem 0.1:}
By Duhamel's principle we have
\begin{equation}
\begin{split}
e^{-t\Delta_g}-e^{-t\Delta_h} 
&=\int_0^t
e^{-s\Delta_g}(\Delta_h-\Delta_g)e^{-(t-s)\Delta_h})\;ds\\
&=\int_0^{t/2}e^{-s\Delta_g}(\Delta_h-\Delta_g)e^{-(t-s)\Delta_h})\;ds\\
&\quad+\int^t_{t/2}e^{-s\Delta_g}(\Delta_h-\Delta_g)e^{-(t-s)\Delta_h})\;ds.
\end{split}
\end{equation}
The integrals converge in the strong operator topology.  By Lemma
\ref{lem6.3} the first integral is a trace class operator.  
In order to prove that the second integral is a trace class operator, it is 
sufficient to prove, that its adjoint with respect to $h$ is of the trace
class.  This adjoint can be written as the strong integral
\begin{equation}\label{6.9}
\int_{\frac{t}{2}}^t \left(e^{-(t-s)\Delta_g}\right)^{\ast_h} \left(
  \Delta_h-(\Delta_g)^{\ast_h}\right)e^{-s\Delta_h} ds. 
\end{equation}
Since $(e^{-(t-s)\Delta_g})^{\ast_h}$ is uniformly bounded in $s$, it follows
again from Lemma \ref{lem6.3} that (\ref{6.9}) is a trace class operator.
\hfill{$\Box$}

\section{Existence and completeness  of wave operators} 
\sectionmark{wave operators}
\setcounter{equation}{0}

In this section we study the wave operators associated to
$(\Delta_g,\Delta_h)$ for equivalent metrics $g$ and $h$.  
\begin{theo}\label{th7.1} Let $g$ and $h$ be two complete metrics of bounded
  curvature on $M$ which satisfy the assumptions of Theorem \ref{th0.1}. 
Let $P_{ac}(\Delta_g)$ be the orthogonal projection onto the absolutely
continuous subspace of $\Delta_g$. 
Then the strong wave operators
\[ W_{\pm}(\Delta_h,\Delta_g) = s-\lim_{t\rightarrow\pm\infty}
e^{it\Delta_h}e^{-it\Delta_g}P_{ac}(\Delta_g) \] 
exist and are complete. In particular, the absolutely continuous parts of
$\Delta_g$ and $\Delta_h$ are unitarily equivalent.   
\end{theo}
\begin{proof} By Theorem \ref{th0.1}, $e^{-t\Delta_g}-e^{-t\Delta_h}$ is
trace class. Then the existence and completeness of the wave operators 
follows from the invariance principle of Birman and Kato
\cite[Chapter X, Theorem 4.7]{Ka}.
\end{proof}

\noindent
{\bf Examples.} We give some examples to demonstrate  Theorem \ref{th0.1}:

\noindent
{\bf 1)} Let $M$ be a manifold with cylindrical ends. Then $\tilde{\imath}$
  is bounded from below, and we may take $b=2$, $a=0$. The condition
  $\beta^{\frac{2}{3}}\in L^1(M)$ is satisfied for $\beta(t) =
  t^{-\frac{3}{2}-\varepsilon}$ for any $\varepsilon>0$.  

\noindent
{\bf 2)} More generally, let $M$ be a manifold with bounded geometry of order
  $2$. (i.e. there is a lower bound for the injectivity radius and the
  covariant derivatives of the curvature of order $\leq 2$ are bounded). Then
  we may choose $x_0 \in M$ arbitrary and let $\beta(t) \leq
  \vol(B_t(x_0))^{-\frac{3}{2}-\varepsilon}$ for any
  $\varepsilon>0$. To see 
  this we first notice that if $M$ is non-compact, the volume of such a
  manifold is infinite. This follows from G\"{u}nthers inequality because we
  may find infinitely many disjoint balls of the same radius. 
Let $a(r):= \frac{\partial}{\partial r}\vol(B_r(x_0))$. Then
$\int_0^1 a(r)\beta(1+r)^{\frac{2}{3}} dr < \infty$ and
\begin{equation*}
\begin{split}
\int_1^\infty a(r)\beta(1+r)^{\frac{2}{3}} dr \leq \int_1^\infty
a(r)\beta(r)^{\frac{2}{3}} dr 
 \leq& \int_1^\infty a(r)\left( \int_0^r a(s) ds
\right)^{-1-\frac{2\varepsilon}{3}} dr \\
&= \int_{\vol(B_1(x_0))}^\infty
t^{-1-\frac{2\varepsilon}{3}} dt < \infty.
\end{split}
\end{equation*}

\noindent
{\bf 3)} Let $M$ be a Riemannian manifold with cusps in the sense of
\cite{Mu1}. Assume that $M$ has bounded curvature.
Then the injectivity radius is
exponentially decreasing in the distance and the volume of $M$ is
finite. Thus we may take $b=1$. It follows $a=1$, and we may take $\beta(t)
= e^{-(\frac{n(n+1)}{2}+4n)ct}$, where $c$ is chosen such that
$\tilde{\imath}(x)\geq Ce^{-cd(x,q)}$. 

\hfill$\square$

The assumptions on $\beta$ in Theorem \ref{th7.1} that guarantee the 
existence of the wave operators are not optimal. Under additional assumptions
on $(M,g)$, the conditions on $\beta$  can be relaxed. For example,
let  $(M,g)$ be a complete manifold which is Euclidean at infinity and
let $h$ be a metric on $M$ which satisfies (\ref{1.12a}), that is
$(M,h)$ is an asymptotically Euclidean manifold.
Then Cotta-Ramusino, Kr\"uger, and Schrader \cite{CKS} proved that the wave
operators $W_{\pm}(\Delta_g,\Delta_h)$ exist. 
The condition (\ref{1.12a}) is weaker then the assumption  which is
necessary in Theorem \ref{th7.1} in this case. The proof is based on Enss's
method \cite{Si}, which applies to this scattering system. An abstract
version of Enss's method has been  developed by Amrein, Pearson 
and Wollenberg \cite{APW}, \cite[16,IV,\S 15]{BW}. This method can
 be applied in  cases where the structure of the continuous spectrum of the
``free Hamiltonian'' is sufficiently well know. To explain this in more detail
we need to introduce some notation. 

Let $C_\infty(\R)$ be the space of all continuous functions on $\R$ that
vanish at infinity. For any closed countable subset $I\subset \R$ let
$C_\infty(\R-I)$ of all functions $f\in C_\infty(\R)$ satisfying $f(x)=0$ for
$x\in I$. A subset $\cA_I$ of the space $C(\R)$ of all bounded continuous 
functions on $\R$ is called multiplicative generating for  $C_\infty(\R-I)$,
if the linear span of the set 
$$\{f\mid f=hg,\;h\in\cA_I,\; g\in C^\infty_c(\R-I)\}$$
 is dense in $C_\infty(\R-I)$ with respect to the norm
$\parallel f\parallel=\sup_{x\in\R}|f(x)|$. The main result of \cite{APW}
 can be stated as follows.

\begin{theo}\label{th7.2}
Let $H$ and $H_0$ be two self-adjoint operators in a Hilbert space $\H$. Let
$R_H(\lambda)$ and $R_{H_0}(\lambda)$ denote the resolvents of $H$ and $H_0$,
respectively. Assume that there exist self-adjoint operators $P_+$ and $P_-$
in $\H$ and a set $\cA_I$ of multiplicative generating functions with respect
to some closed countable subset $I\subset \R$ satisfying the following
properties 
\begin{enumerate}
\item $P_{\mathrm{ac}}(H_0)=P_++P_-$ and
  $\slimit_{t\to\pm\infty}e^{itH_0}P_\mp e^{-itH_0}P_{\mathrm{ac}}(H_0)=0.$
\item $(\Id-P_{\mathrm{ac}}(H_0))\alpha(H_0)$ is compact for all $\alpha\in
\cA_I$.
\item $R_{H}(i)-R_{H_0}(i)$ is compact.
\item $\int_0^{\pm\infty}\parallel(R_{H}(i)-R_{H_0}(i))e^{-itH_0}
\alpha(H_0)P_\pm\parallel \;dt<\infty$ for all $\alpha\in\cA_I$. 
\end{enumerate}
Then the wave operators $W_\pm(H,H_0)$ exist and are complete. Moreover $H$
and $H_0$ have no singularly continuous spectrum and each eigenvalue of $H$
and $H_0$ in $\R-I$ is of finite multiplicity. These eigenvalues accumulate at
most at points of $I\cup\{\pm\infty\}$.
\end{theo}

For the proof see Corollary 19 in \cite[16,IV,\S 15]{BW}.

As example, we consider a manifold $X$ with  cusps as defined in \cite{Mu1}. 
For simplicity we assume that $X$ has a single cusp. Then  $X$ is a complete
Riemannian manifold of dimension $n+1$ that admits a decomposition 
$$X=M\cup_Y Z$$
in a compact Riemannian  manifold $M$ with boundary $Y$ and a half-cylinder
$Z=[1,\infty)\times Y$, and $M$ and $Z$ are glued along their common boundary 
$Y$. The metric $g$ on $X$ is such that its restriction to $Z$ is given by
\begin{equation}\label{7.0}
g^Z=u^{-2}(du^2+g^Y),
\end{equation}
where $g^Y$ denotes the metric of $Y$. The metric $g$ is the fixed background
metric and we consider perturbations $h$ of $g$.  
 As free Hamiltonian $H_0$  we are taking a modification
of the Laplacian $\Delta_g$ which is defined as follows.  
We regard $Y$ as a hypersurface in $X$
that separates $X$ into
 $M$ and $Z$. Let $C^\infty_0(X-Y)$ be the subspace of
all $f\in C^\infty_c(X)$ that vanish in a neighborhood of $Y$. Let
$\Delta_0$ denote Friedrichs's extension of 
$$\Delta_g \colon C^\infty_0(X-Y)\to L^2(X).$$
To begin with we need to study the spectrum of
$\Delta_0$. 
With respect to the decomposition $L^2(X)=L^2(M)\oplus L^2(Z)$ we have 
\begin{equation}\label{7.1}
\Delta_0=\Delta_{M,0}\oplus\Delta_{Z,0},
\end{equation}
where $\Delta_{M,0}$ and $\Delta_{Z,0}$ are the Dirichlet Laplacians
on $M$ and $Z$, respectively. 
Since $M$ is compact, $\Delta_{M,0}$ has pure point spectrum. Let 
$$L^2_0(Z):=\bigl\{f\in L^2(Z)\colon \int_Y f(u,y)\;dy=0\;\mathrm{for\; almost
  \;all}\; u\in[1,\infty)\bigr\}.$$
The orthogonal complement $L_0^2(Z)^\perp$ of $L^2_0(Z)$ in $L^2(Z)$ consists 
of functions which are independent of $y\in Y$ and therefore, can be
identified with $L^2([1,\infty),u^{-(n+1)}du)$. The decomposition
\begin{equation}\label{7.2}
L^2(Z)=L^2_0(Z)\oplus L^2_0(Z)^\perp
\end{equation}
is invariant under $\Delta_{Z,0}$. 
\begin{lem}\label{l7.3}
The restriction of $\Delta_{Z,0}$ to $L^2_0(Z)$ has a compact resolvent. In
particular, $\Delta_{Z,0}$ has pure point spectrum.
\end{lem}
\begin{proof}
Let $\Delta_Y$ be the Laplacian of $Y$. Let $\{\phi_j\}_{j=0}^\infty$ be an
orthonormal basis of eigenfunctions of $\Delta_Y$ with eigenvalues $0=\lambda_0
<\lambda_1\le\lambda_2\le\cdots$. Let $f\in C^\infty_c(Z)\cap L^2_0(Z)$.
Then $f$ has an expansion of the form
$$f(u,y)=\sum_{k=1}^\infty a_k(u)\phi_k(y),$$
where the series  converges in the $C^\infty$-topology. Let $b>1$ and put
$Z_b=[b,\infty)\times Y$. Let $C=\lambda_1^{-1}$. Then we have
\begin{equation}\label{7.3}
\parallel f\parallel^2_{L^2(Z_b)}=\sum_{k=1}^\infty\int_b^\infty
|a_k(u)|^2\, \frac{du}{u^{n+1}}\le \frac{C}{b^2}\sum_{k=1}^\infty \lambda_k
\int_b^\infty |a_k(u)|^2\, \frac{du}{u^{n-1}}.
\end{equation}
Now observe that the Laplacian $\Delta_Z$ with respect to the metric (\ref{7.0})
equals
\begin{equation}\label{7.3b}
-u^2\frac{\partial^2}{\partial u^2}+nu\frac{\partial}{\partial u}
+u^2\Delta_Y.
\end{equation}
Moreover, since $a_k\in C^\infty_c((1,\infty))$, we have
$$\int_1^\infty\left(-u^2 a_k^{\prime\prime}(u)+
nu a_k^\prime(u)\right)\overline{a_k(u)}\;\frac{du}{u^{n+1}}=\int_1^\infty|a_k^\prime(u)|^2u^{1-n}\;du\ge 0.$$
This together with (\ref{7.3}) implies
\begin{equation}\label{7.3a}
\parallel f\parallel^2_{L^2(Z_b)}\le\frac{C}{b^2}\langle\Delta_Zf,
f\rangle_{L^2(Z)}= \frac{C}{b^2}\parallel \nabla
f\parallel^2_{L^2(Z)}\le \frac{C}{b^2}\parallel f\parallel_{H^1(Z)}^2.
\end{equation}
Let $H^1_0(Z):= H^1(Z)\cap L^2_0(Z)$.
By continuity,  (\ref{7.3})  holds for all $f\in H^1_0(Z)$. By Rellich's 
lemma, the embedding
$$i_b\colon H^1(Z-Z_b)\cap L^2_0(Z-Z_b)\to L^2(Z)$$
is compact. It follows from (\ref{7.3a}) that as $b\to\infty$, $i_b$ converges
strongly to the embedding
$$i\colon H^1_0(Z)\to L^2(Z).$$
Hence $i$ is compact which implies the lemma.
\end{proof}
Let 
$$D_0:= -u^2\frac{d^2}{du^2}+nu\frac{d}{du}\colon C_c^\infty((1,\infty))\to 
L^2([1,\infty),u^{-(n+1)}du)$$
and let $L_0$ be the self-adjoint extension of 
$D_0$ with respect to Dirichlet boundary conditions at $1$.
By (\ref{7.3b}), the restriction of $\Delta_{Z,0}$ to  $L^2_0(Z)^\perp\cong
L^2([1,\infty),u^{-(n+1)}du)$ is equivalent to $L_0$.  The spectrum of $L_0$ 
is absolutely continuous and equals $[n^2/4,\infty)$. Thus we get the
following lemma.
\begin{lem}\label{l7.4}
The spectrum of $\Delta_0$ is the union of a pure point point spectrum and
an absolutely continuous spectrum. The point spectrum consists of eigenvalues
of finite multiplicity $0<\lambda_1<\lambda_2<\cdots\to\infty$.
The absolutely continuous spectrum is equal to $[n^2/4,\infty)$ and the
absolutely continuous  part $\Delta_{0,\ac}$ of $\Delta_0$ is equivalent 
to $L_0$. 
\end{lem}
Let $\varepsilon>0$ and let $\beta(t)=e^{-\varepsilon t}$. Let $h$ be a
complete  
metric on $X$. We put 
$$H:=\Delta_h\quad\mathrm{and}\quad H_0:=\Delta_0.$$ 
Since $H$ and $H_0$ are positive operators, we can replace $i$ by $-1$ in
 Theorem \ref{th7.2}. So let
\begin{equation}\label{7.4}
R_g:=(\Delta_g+\Id)^{-1},\;\;R_h:=(\Delta_h+\Id)^{-1}\,\,\mathrm{and}\;\;
R_0:=(\Delta_0+\Id)^{-1}. 
\end{equation}
First we have the following lemma.
\begin{lem}\label{l7.5}
Suppose that $h\sim_\beta^2 g$. Then $R_h-R_0$ is a compact
operator. 
\end{lem}
\begin{proof}  Since $Y\subset X$ is a compact
hypersurface, it follows that $R_g-R_0$ is a compact operator. So it suffices
to show that $R_h-R_g$ is compact. We have
\begin{equation}\label{7.5}
R_h-R_g=-R_g(\Delta_h-\Delta_g)R_h.
\end{equation}
By Lemma \ref{l5.2} we have
\begin{equation}\label{7.6}
\Delta_h-\Delta_g=\sum_{j=0}^2\xi_j\circ(\nabla^h)^j
\end{equation}
and $\xi_j$ satisfies
\begin{equation}\label{7.7}
|\xi_j(x)|\le C e^{-\varepsilon d(x,x_0)},\quad x\in X.
\end{equation}
Now $R_h\colon L^2(X)\to W^2(X)$ is continuous. Therefore by (\ref{7.6})
and (\ref{7.7}) it follows that 
$$(\Delta_h-\Delta_g)R_h\colon L^2(X)\to L^2(X)$$
is a bounded operator. Using again that $R_g-R_0$ is compact, it
follows from (\ref{7.5}) that it suffices to show that
$R_0(\Delta_h-\Delta_g)R_h$ is a compact operator. 

For $a>1$  let 
$$X_a=M\cup_Y([1,a]\times Y).$$
Denote by $\chi_a$ the characteristic function of $X_a$ in $X$.
We claim that $R_0\chi_a$ is a compact operator. By (\ref{7.1}) 
we have
$$R_0=(\Delta_{M,0}+1)^{-1}\oplus(\Delta_{Z,0}+1)^{-1}.$$
Since $M$ is compact, $(\Delta_{M,0}+1)^{-1}$ is compact. Let $\chi_{[1,a]}$
be the characteristic function of the interval $[1,a]$ in $[1,\infty)$.
By Lemma \ref{7.3} it suffices to show that
$(L_0+\Id)^{-1}\chi_{[1,a]}$ is compact as operator in
$L^2([1,\infty),u^{-(n+1)} 
du)$. The kernel
$g(u,u^\prime)$ of $(L_0+\Id)^{-1}$ is given by
\begin{equation}\label{7.8}
g(u,u^\prime)=\frac{(uu^\prime)^{n/2}}{\sqrt{n^2/4+1}}
\begin{cases}(u^\prime/u)^{(n^2/4+1)^{1/2}},\,u>u^\prime;\\
(u/u^\prime)^{(n^2/4+1)^{1/2}},\,u^\prime>u.\end{cases}
\end{equation}
From this formula follows that
$g(u,u^\prime)$ is bounded on 
$[1,\infty)\times [1,a]$, and therefore square integrable with respect to the
measure $u^{-(n+1)}du$. This implies that $(L_0+\Id)^{-1}\chi_{[1,a]}$ is
a compact operator and hence, $R_0\chi_a$ is compact for all $a>1$. 

Let $M_{(1-\chi_a)\beta}$ denote the multiplication operator by
$(1-\chi_a)\beta$. Using
(\ref{7.6}) and (\ref{7.7}), we get
\begin{equation}\label{7.9}
\begin{split}
\parallel R_0&(1-\chi_a)(\Delta_h-\Delta_g)R_h\parallel\\
&\le
C\left(\sum_{j=0}^2\parallel(\nabla^h)^jR_h\parallel 
\right)\cdot \parallel R_0\parallel\cdot\parallel M_{(1-\chi_a)\beta}
\parallel.
\end{split}
\end{equation}
Let $Z_a=[a,\infty)\times Y$. Then
$$\parallel M_{(1-\chi_a)\beta}\parallel\le \sup_{x\in Z_a}\beta(x)=
\sup_{x\in Z_a}e^{-\varepsilon d(x,x_0)}.$$
Now observe that there exists $C_1>0$ such that for all $(u,y)\in Z_a$ we have
$$d((u,y),x_0)\ge d((u,y),(1,y))-C=\log u -C_1.$$
Hence together with (\ref{7.9}) we get
$$\parallel R_0(1-\chi_a)(\Delta_h-\Delta_g)R_h\parallel\le C_2
a^{-\varepsilon}.$$
Thus $R_0(\Delta_h-\Delta_g)R_h$ can
be approximated in the operator norm by compact operators and hence, is a
compact operator.
\end{proof}

Next we construct self-adjoint projections $P_\pm$ which satisfy the 
conditions of Theorem \ref{th7.2}. Let
\begin{equation}\label{7.10}
e(u,\lambda):= u^{n/2+i\lambda}-u^{n/2-i\lambda},\quad u\in[1,\infty),\;
\lambda\in\R.
\end{equation}
Then $e(u,\lambda)$ satisfies
$$D_0 e(u,\lambda)=(n^2/4+\lambda^2)e(u,\lambda),\quad e(1,\lambda)=0.$$
Thus $e(u,\lambda)$ is the generalized eigenfunction for $L_0$. For
$\varphi\in C^\infty_c(1,\infty))$ set
$$\hat\varphi(\lambda):=\frac{1}{2\pi}\int_1^\infty e(u,\lambda)\varphi(u)\;
\frac{du}{u^{n+1}}.$$
The map $\varphi\mapsto \hat\varphi$ extends to an isometry
$$F\colon L^2([1,\infty),u^{-(n+1)}du)\to L^2(\R^+)$$
such that
$$F\circ L_0\circ F^*=\widetilde L_0,$$
where $\widetilde L_0$ is the multiplication operator by $(n^2/4+\lambda^2)$.
Let 
$$U\colon L^2(\R^+)\to L^2([n^2/4,\infty))$$
be defined by
$$(Uf)(\lambda)=\frac{f(\sqrt{\lambda-n^2/4})}
{\sqrt{2}(\lambda-n^2/4)^{1/4}}.$$ 
Then $U$ is an isometry such that $U\circ \widetilde L_0\circ
U^*=\widehat L_0$, where $\widehat L_0$ is the multiplication operator by
$\lambda$. Thus $U\circ F$ provides the spectral resolution of
$L_0=\Delta_{0,\ac}$. Let
$$J\colon L^2([n^2/4,\infty))\to L^2(\R)$$ 
denote the inclusion, let $\cF\colon L^2(\R)\to L^2(\R)$ be the Fourier
transform, and let $\chi_\pm$ denote the characteristic function of
$[0,\infty)$ and $(-\infty,0]$, respectively. Set
$$\widetilde P_\pm:= J^*\cF \chi_\pm \cF^*J.$$
Then $\widetilde P_++\widetilde P_-$ is the identity of $L^2([n^2/4,\infty))$.
Let $A=-id/du$, regarded as self-adjoint operator in $L^2(\R)$. Then
$$\widetilde P_\pm e^{-it\widehat L_0}=J^*\cF \chi_\pm e^{-itA} \cF^*J.$$
Let $f\in L^2(\R)$. Using the Fourier transformation, it follows that
$(e^{-itA}f)(u)=f(u-t)$. Thus we get
$$\parallel \chi_\pm e^{-itA}f\parallel^2=\pm\int_{-t}^{\pm\infty}|f(u)|^2\,du
\to 0$$
as $t\to \mp\infty$. Hence we get
\begin{equation}\label{7.11}
\slimit_{t\to\pm\infty}e^{it\widehat L_0}\widetilde P_\mp e^{-it\widehat L_0}
=0.
\end{equation}
Now put 
$$P_\pm:=F^*U^*\widetilde P_\pm UF$$
on $L^2([1,\infty),u^{-(n+1)}du)$ and set $P_\pm:=0$ on the orthogonal
complement of $L^2_0(Z)^\perp=L^2([1,\infty),u^{-(n+1)}du)$ in $L^2(X)$.
Then $P_\pm$ are self-adjoint projections that satisfy 
$$P_+ +P_-=P_\ac(\Delta_0).$$
Furthermore we have
$$e^{itH_0}P_\pm e^{-it\Delta_0}P_\ac(\Delta_0)=F^*U^*e^{it\widehat
  L_0}\widetilde P_\pm e^{-it\widehat L_0}UF.$$
So it follows from (\ref{7.11}) that
$$\slimit_{t\to\pm\infty}e^{itH_0}P_\mp e^{-itH_0}P_\ac(H_0)=0.$$
Thus condition (1) of Theorem \ref{th7.2} is satisfied. Let $I=\{n^2/4\}$.
and put
$$\cA_I:=C^\infty_c(\R-I).$$
Then it is clear that $\cA_I$ is multiplicative generating for
$C_\infty(\R-I)$. By  Lemma \ref{l7.4}, $\Delta_0$ has
pure point spectrum in the subspace $(\Id-P_\ac(\Delta_0))L^2(X)$ consisting
of eigenvalues of finite multiplicity with no finite points of accumulation.
Let $\alpha\in \cA_I$. Then  $(\Id-P_\ac(\Delta_0))\alpha(\Delta_0)$
is a finite rank operator. This is condition (2) of Theorem \ref{th7.2}.
Condition (3) holds by Lemma \ref{l7.5}. It remains to verify condition
(4). 

Given $t>0$, let $\chi_t$ be the characteristic function of
$[e^t,\infty)\times Y$ in $X$. Let $\delta>0$. We have
\begin{equation}\label{7.12}
\begin{split}
\parallel(R_h&-R_{0})e^{it\Delta_0}\alpha(\Delta_0)
P_\pm\parallel\\
&\le\parallel R_h-R_{0}\parallel\cdot \parallel 
(1-\chi_{\delta t})
 e^{it\Delta_0}\alpha(\Delta_0)P_\pm\parallel \\
&\quad+ \parallel (R_h-R_{0})\chi_{\delta t}\parallel\cdot
\parallel \alpha(\Delta_0)\parallel.
\end{split}
\end{equation}
We will prove that for each $\alpha\in C_c^\infty(\R-\{n^2/4\})$ there exists
$\delta>0$ such that the right hand side is an integrable function of
$t\in\R^+$. 
To estimate the first term on the right hand side we need the following 
auxiliary result.
\begin{lem}\label{l7.6} Let $a\in\R$ and let $f\in C^\infty_c(\R-\{a\})$.
Let $\varepsilon>0$ such that $f(\lambda^2+a)=0$ for $|\lambda|<\varepsilon$.
 Then for every $m\in\N$  there exists $C>0$ such 
that for $t\in\R-\{0\}$ and $|u|<\varepsilon|t|/2$ one has
$$
\bigg|\int_0^\infty e^{2iu\lambda+it\lambda^2}f(\lambda^2+a)\,d\lambda\bigg|
\le C|t|^{-m}.
$$
\end{lem}
\begin{proof}
Let $t\not=0$ and set $x=u/t$. Then the left hand side of the inequality
equals
\begin{equation*}
\begin{split}
\bigg|\int_0^\infty e^{it(\lambda+x)^2}&f(\lambda^2+a)\,d\lambda\bigg|\\
&=(2t)^{-m}\bigg|\int_0^\infty e^{it(\lambda+x)^2}
\left (\frac{1}{\lambda+x}\frac{d}{d\lambda}-\frac{1}{(\lambda+x)^2}\right)^m 
f(\lambda^2+a)\,d\lambda\bigg|.
\end{split}
\end{equation*}
Now assume that $|u|<\varepsilon |t|/2$. Then $|x|<\varepsilon/2$. On the
other hand, we have $f(\lambda^2+a)=0$ for $|\lambda|<\varepsilon$. Thus if
$f(\lambda^2+a)\not=0$, then we have
$|\lambda+x|\ge|\lambda|-|x|>\varepsilon/2$. Hence the right hand side can be
estimated by $C|t|^{-m}$. 
\end{proof}

Let $\varphi\in L^2([1,\infty),u^{-(n+1)}du)=P_\ac(\Delta_0)(L^2(X))$. Then
\begin{equation}\label{7.13}
\begin{split}
\bigl(e^{-it\Delta_0}&\alpha(\Delta_0)\varphi\bigr)(u)\\
&=\frac{1}{2\pi}\int_0^\infty e(u,n/2-i\lambda)e^{-it(\lambda^2+n^2/4)}
\alpha(\lambda^2+n^2/4)(F\varphi)(\lambda)\,d\lambda.
\end{split}
\end{equation}
Let  $v\in C^\infty_c((1,\infty))$. Put $\varphi=P_+v$ and
$w=\cF^*JUFv$. Then $w\in L^1(\R)$ and $F\varphi=U^*J^*\cF(\chi_+w)$. Using
the definition of $U$, $J$ and $\cF$, we get
$$(F\varphi)(\lambda)=\sqrt{2\lambda}\int_0^\infty e^{-is(\lambda^2+n^2/4)} 
w(s)\,ds.$$
Assume that $t>0$. If we insert this
expression into the right hand side of  (\ref{7.13}) and  switch the order
of integration, we obtain
\begin{equation}\label{7.14}
\begin{split}
\bigl(&e^{-it\Delta_0}\alpha(\Delta_0)P_+v\bigr)(u)\\
&=\frac{1}{\sqrt{2}\pi}\int_0^\infty w(s)\int_0^\infty e(u,n/2-i\lambda)
e^{-i(t+s)(\lambda^2+n^2/4)}\alpha(\lambda^2+n^2/4)\sqrt{\lambda}
\,d\lambda\,ds. 
\end{split}
\end{equation}
Now there exists $\varepsilon>0$ such that $\alpha(\lambda^2+n^2/4)=0$ for
$|\lambda|<\varepsilon$. Assume that $|\log(u)|<\varepsilon t/2$. Using the 
definition (\ref{7.10}) of $e(u,\lambda)$ and Lemma \ref{l7.6}, it follows
that there exists $C>0$ such that
\begin{equation}\label{7.15}
\begin{split}
\bigg|e^{-it\Delta_0}&\alpha(\Delta_0)P_+v(u)\bigg|^2
\le C \parallel w\parallel^2  u^n t^{-3}\le C  \parallel v\parallel^2u^n
t^{-3}. 
\end{split}
\end{equation}
Thus for every $\alpha\in C_c^\infty(\R-\{n^2/4\})$ there exist $C>0$ and
$\delta>0$ such that for $t>\delta^{-1}$ one has
$$\parallel (1-\chi_{\delta t})e^{-it\Delta_0}\alpha(\Delta_0)P_+\parallel
\le C \,t^{-3}\int_1^{e^{\delta t}}\frac{du}{u} =C\delta t^{-2}.$$
Similarly one can show that
$$\parallel (1-\chi_{\delta t})e^{-it\Delta_0}\alpha(\Delta_0)P_-\parallel
\le C \,t^{-3}\int_1^{e^{\delta t}}\frac{du}{u} =C\delta t^{-2},\quad
t>\delta^{-1}.$$
Hence for this choice of $\delta$, the first term on the right hand side of
(\ref{7.12}) is an integrable function of $t\in\R^+$.

Now consider the second term on the right hand side of (\ref{7.12}). We have
\begin{equation}\label{7.16}
\begin{split}
\parallel (R_h-R_{0})\chi_{\delta t}\parallel
\le \parallel (R_h-R_{g})\chi_{\delta t}\parallel+
\parallel (R_g-R_{\Delta_0})\chi_{\delta t}\parallel.
\end{split}
\end{equation}
Let $M_{\chi_{\delta t}\beta}$ denote the multiplication operator by
$\chi_{\delta t}\beta$.  By (\ref{7.5}) - (\ref{7.7}) we get
\begin{equation}\label{7.17}
\begin{split}
\parallel (R_h-R_{g})\chi_{\delta t}\parallel&\le \parallel R_g
\parallel\cdot \parallel\chi_{\delta t}(\Delta_h-\Delta_g)R_h\parallel\\ 
&\le C \parallel M_{\chi_{\delta t}\beta}\parallel\left(\sum_{j=0}^2\parallel
(\nabla^h)^jR_h\parallel\right)\le C_1 e^{-\varepsilon\delta t}. 
\end{split}
\end{equation}
It remains to estimate the second term on the right of (\ref{7.16}).
Let $\psi\in C^\infty(\R)$ such that $f(u)=0$, if $u\le 2$, and $f(u)=1$, if 
$u\ge 3$. Define $f\in C^\infty(Z)$ by $f(u,y)=\psi(u)$ and extend $f$ by
zero to a smooth function on $X$. Then we have
\begin{equation*}
\begin{split}
R_g-R_{0}&=(f-1)R_{0}- R_g((\Delta_g+\Id)(f R_{0})-\Id).
\end{split}
\end{equation*}
Observe that
$$(\Delta_g+\Id)(fR_{0})-\Id=
f-1+2\nabla f\cdot\nabla R_{0}+\Delta f\cdot R_{0}.$$ 
Moreover note that $(f-1)\chi_{\delta t}=0$ if $t\gg0$. Thus
\begin{equation}\label{7.18}
\begin{split}
(R_g-R_{0})\cdot\chi_{\delta t}
=(f-1)\cdot R_{0}\cdot\chi_{\delta t}
- R_g(2\nabla f\cdot\nabla R_{0}\cdot\chi_{\delta t}
+\Delta f\cdot R_{0}\cdot\chi_{\delta t}) 
\end{split}
\end{equation}
for $t\gg0$. It follows from (\ref{7.1}) that $R_{0}\cdot\chi_{\delta t}$
acts in 
$L^2(Z)$ and preserves the decomposition (\ref{7.2}). Moreover $\parallel
R_{0}\cdot\chi_{\delta t}|_{L^2_0(Z)}\parallel=
\parallel \chi_{\delta t}\cdot R_{0}|_{L^2_0(Z)}\parallel$.
Let $\varphi\in L^2_0(Z)$. Then $R_{0}\varphi\in
L^2_0(Z)\cap H^2(Z)$ and by (\ref{7.3}) we obtain
\begin{equation}\label{7.19}
\begin{split} 
\parallel \chi_{\delta t}R_{0}\varphi\parallel\le C e^{-2\delta t}
\parallel R_{0}\varphi\parallel_1
\le C e^{-2\delta t}\parallel\varphi\parallel.
\end{split}
\end{equation}
On the orthogonal complement $L^2_0(Z)^\perp$, the kernel of $R_{0}$
 is given by (\ref{7.8}). Let $h\in C^\infty_c(Z)$. Then it follows from 
(\ref{7.8}) that
\begin{equation}\label{7.20}
\begin{split}
\parallel h\cdot R_{0}\cdot\chi_{\delta t}
|_{L^2_0(Z)^\perp}\parallel\le C e^{-\delta t\sqrt{n^2/4+1}}
\le C e^{-\delta t}.
\end{split}
\end{equation}
Combining (\ref{7.19}) and (\ref{7.20}) we obtain
$$\parallel h \cdot R_{0}\cdot\chi_{\delta t}\parallel\le C e^{-\delta t}.$$
Similar estimations hold for $\nabla R_{0}$. This proves that
the second term on the right hand side of (\ref{7.12}) is an integrable 
function of $t\in\R^+$. This is condition (4) of Theorem \ref{th7.2}.
Summarizing we have proved the following theorem.

\begin{theo}\label{th7.7}
Let $(X,g)$ be a manifold with cusps and let $\Delta_0$ be 
defined by (\ref{7.1}). Let $\varepsilon>0$ and put $\beta(u)=e^{-\varepsilon
  u}$, $u\in\R$. Let $h$ be a complete metric on $X$ such that
$h\sim_\beta^2 g$. Then we have
\begin{enumerate}
\item The wave operators $W_\pm(\Delta_h,\Delta_0)$ exist and are complete.
\item $\Delta_h$ has no singularly continuous spectrum.
\end{enumerate}
\end{theo}

\begin{cor}\label{c7.8}
Let $g$ and $h$ be as above. Then the wave operators
$W_\pm(\Delta_h,\Delta_g)$ exist and are complete.
\end{cor}
This is a considerable improvement of the result that we get from Theorem 
\ref{th7.1} in this case.

\smallskip
{\bf Remark.} Other cases of complete manifolds $(M,g)$ with a sufficiently
explicit structure at infinity can be treated in the same way. This includes,
for example, manifolds with cylindrical ends and asymptotically Euclidean
manifolds.

\sectionmark{$beta$-Equiv. and Analytic Cont. of the Resolvent.}
\section{$\beta$-Equivalence and Analytic Continuations of the Resolvent} 
\sectionmark{$beta$-Equiv. and Analytic Cont. of the Resolvent.}
\setcounter{equation}{0}

In this section we study the  existence of an analytic continuation
 of the resolvent in weighted $L^2$-spaces. 
Provided that such a continuation exists, we are able to study the behavior 
of the absolutely continuous spectrum under perturbation in more detail. 
 The method is a modification of the method used in \cite{Mu2}. 

\begin{defn} Let ${\cB}$ be a Banach space, $\Omega\subset \C$ a domain and
 $F:\Omega \mapsto {\cB}$ a meromorphic function. Let $\Sigma$ be a Riemann
surface and let $\pi:\Sigma \to \C$ be a ramified covering. A meromorphic 
continuation of $F$ to $\Sigma$ is a meromorphic  function 
$\tilde{F}:\Sigma \to {\cB}$ such that
\begin{itemize}
\item[a)] There exists $\tilde{\Omega} \subseteq \Sigma$ such that $\pi\colon
\tilde{\Omega}\to \Omega$ is biholomorphic.
\item[b)] $F\circ \pi = \tilde{F}$ on $\tilde{\Omega}$.
\end{itemize}
\end{defn}

\begin{defn} Let $\delta$ be a function of moderate decay and let 
$p\in\N$. By $H^{-p}_{\delta^{-1}}$ we denote the dual space of $H^p_\delta$, 
with respect to  the extension of the $L^2$-pairing.
\end{defn}

\begin{lem} \label{l8.3}  Let $\zeta(u)$  be a non-increasing 
continuous function on $[1,\infty)$  with $\zeta(u)\to 0$ as $u\to\infty$ and
 let
$\delta$ be a weight function. 
Then  the canonical inclusion
$j\colon L^2_{\delta\zeta^{-1}}(M) \to H^{-2}_\delta(M)$ is compact.
\end{lem}
\begin{proof} It is enough to prove, that the adjoint
$\jmath^\ast:H^2_{\delta^{-1}}(M) \to L^2_{\delta^{-1}\zeta}(M)$ is compact. 
For $k\in \N$ let
$$\Omega_k=\{x\in M\mid \zeta(1+d(x,x_0))\ge 1/k\}.$$
Then each $\Omega_k$ is a compact subset of $M$. Let $P_k$ be the 
multiplication operator by the characteristic function of $\Omega_k$. 
By Rellich's lemma, $\jmath^\ast P_k$ is compact. For $f\in
H^2_{\delta^{-1}}(M)$ 
we have
$$\int_{M-\Omega_k}|f(x)|^2\delta^{-1}(x)\zeta(x)\;dx\le \frac{1}{k}\parallel f
\parallel^2_{H^2_{\delta^{-1}}}.$$
Thus $\jmath^\ast P_k$ converges to $\jmath^\ast$ in the operator topology.
Hence $\jmath^\ast$ is compact. 
\end{proof}

Let $\delta,\rho$ be functions of moderate decay. Then 
$L^2_{\delta^{-1}}(M)\subset L^2(M)$ and $H^2(M)\subset H^2_\rho(M)$. Thus for 
$\lambda\in\C-[0,\infty)$, the resolvent $(\Delta-\lambda)^{-1} \colon L^2(M)
\to H^2(M)$ may be 
regarded as a bounded operator
$$(\Delta-\lambda)^{-1}\colon  L^2_{\delta^{-1}}(M)\to H^2_\rho(M).$$  
Denote by $\mathcal{L}(L^2_{\delta^{-1}}(M),H^2_\rho(M))$ the Banach space of
all 
bounded operators from $L^2_{\delta^{-1}}(M)$ into $H^2_\rho(M)$, equipped with
the strong operator norm.

\begin{theo}\label{th8.4}
Let $g,h$ be  complete Riemannian metrics on $M$ with bounded
curvature of order 2. Let $\beta$, $\delta$, $\zeta$ and $\rho$ 
be functions of moderate decay
on $M$ such that
\begin{equation}\label{8.1} 
\beta^2(x) \leq C\tilde{\imath}^{4n}_g(x)\rho(x)\delta(x)\zeta(x),
 \quad x\in M, 
\end{equation}
and $g\sim^2_\beta h$.
Let $\Omega\subset \C-[0,\infty)$ be open. Assume that there is a 
Riemann surface $\Sigma$ and a covering $\Sigma\to\Omega$  such that
the operator valued function
$$\lambda\in\Omega\mapsto (\Delta_g-\lambda)^{-1}\in 
\mathcal{L}(L^2_{\delta^{-1}}(M,g),H^2_\rho(M,g))$$
admits an analytic continuation to a meromorphic function 
$$\lambda\in\Sigma\to R_g(\lambda) 
\in\mathcal{L}(L^2_{\delta^{-1}}(M,g),H^2_\rho(M,g))$$
with finite rank residues. Then 
$$\lambda\in\Omega\mapsto (\Delta_h-\lambda)^{-1}\in 
\mathcal{L}(L^2_{\delta^{-1}}(M,h),H^2_\rho(M,h))$$ 
also admits a meromorphic continuation  to $\Sigma$ with finite rank residues.
\end{theo}
\begin{proof} By assumption, $\beta\tilde{\imath}^{-2n}$ is bounded.
Hence by Corollary \ref{c5.3}, $H^2(M,g)$ and 
$H^2(M,h)$ are equivalent and therefore, by duality, $H^{-2}(M,g)$ and
 $H^{-2}(M,h)$ are also equivalent. Let $\lambda\in \C-[0,\infty)$. Then 
$$K(\lambda):=(\Delta_g-\lambda)^{-1}(\Delta_h-\Delta_g)$$
is a bounded operator in $L^2(M)$. Moreover 
$\Id+K(\lambda)=(\Delta_g-\lambda)^{-1}(\Delta_h-\lambda)$
has a bounded inverse in $L^2(M)$ which is given by
$$(\Id+K(\lambda))^{-1}=(\Delta_h-\lambda)^{-1}(\Delta_g-\lambda).$$
Thus for $\lambda\in \C-[0,\infty)$ we have
\begin{equation}\label{8.2}
(\Delta_h-\lambda)^{-1}=(\Id+K(\lambda))^{-1}(\Delta_g-\lambda)^{-1}.
\end{equation}
By Corollary \ref{c4.3} there exists $\lambda\in\C-[0,\infty)$ such that
$(\Delta_h-\lambda)^{-1}$ extends to a bounded operator in $L^2_\rho(M)$. 
By Lemma \ref{l4.5} it follows that $(\Delta_h-\lambda)^{-1}$ maps
$L^2_\rho(M)$ 
into $H^2_\rho(M)$. Moreover by definition $\Delta_g-\lambda$ is a bounded 
operator
of $H^2_\rho(M)$ to $L^2_\rho(M)$. Hence $(\Id+K(\lambda))^{-1}$ extends to a
bounded operator in $H^2_\rho(M)$. Let $\mu\in \Omega$. Then
\begin{equation}
\begin{split}
 \Id+K(\mu) &= \left(\Id+K(\lambda)\right) - \left\{
   \left(\Id+K(\lambda)\right)- 
\left(\Id+K(\mu)\right)\right\} \\
&= \left(\Id+K(\lambda)\right) - \left\{ K(\lambda)-K(\mu) \right\}\\ 
 &= \left(\Id+K(\lambda)\right) - (\lambda-\mu)(\Delta_g-\mu)^{-1}
(\Delta_g-\lambda)^{-1}(\Delta_h-\Delta_g). 
\end{split}
\end{equation}

By Corollary \ref{c4.3} we may choose $\lambda$ such that 
$(\Delta_g-\lambda)^{-1}$ extends to a bounded operator in 
$L^2_{\delta\zeta}(M)$. By duality, and Lemma \ref{l4.5}, it defines 
a bounded operator 
$$(\Delta_g-\lambda)^{-1}\colon L^2_{\delta^{-1}\zeta^{-1}}(M)\to 
H^2_{\delta^{-1}\zeta^{-1}}(M).$$  
Using Lemma \ref{l3.6}, Lemma \ref{l5.2} and the assumption on $\beta$, it
follows that the operator $(\Delta_g-\lambda)^{-1}(\Delta_h-\Delta_g)$ is the
composition of
the following chain of bounded operators
\begin{equation}
\begin{split}
H^2_\rho(M)\to W^2_{\tilde{\imath}^{4n}\rho}(M)&\stackrel{\Delta_h-\Delta_g}
{\longrightarrow} L^2_{\beta^{-2}\tilde{\imath}^{4n}\rho}(M)\to\\ 
&L^2_{\delta^{-1}\zeta^{-1}}(M)\stackrel{(\Delta_g-\lambda)^{-1}}{\longrightarrow}
H^2_{\delta^{-1}\zeta^{-1}}(M)\stackrel{j}{\longrightarrow} 
L^2_{\delta^{-1}}(M).
\end{split}
\end{equation}
By Lemma \ref{l8.3}, the inclusion $j$ is a compact. Hence
$$(\Delta_g-\lambda)^{-1}(\Delta_h-\Delta_g)\colon H^2_\rho(M)\to 
L^2_{\delta^{-1}}(M)$$
is compact operator. Set
\begin{equation}
H_\lambda(\mu) = (\lambda-\mu)R_g(\mu)\circ(\Delta_g-\lambda)^{-1}
(\Delta_h-\Delta_g),\quad \mu\in\Sigma. 
\end{equation}
Then $H_\lambda(\mu)$, $\mu\in\Sigma$,  is a meromorphic family of compact 
operators and
\begin{equation} \Id+K(\mu) = \left(\Id+K(\lambda)\right)\left
\{ \Id - \left(\Id+K(\lambda)\right)^{-1}H_\lambda(\mu)\right\}. 
\end{equation}
It then follows from \cite{St}, that $(\Id+K(\mu))^{-1}$ exists except 
for on a discrete set and is meromorphic in $\mu$. Thus, we may define
\begin{equation} 
R_h(\mu)= (\Id+K(\mu))^{-1}\circ R_g(\mu). 
\end{equation}
By (\ref{8.2}) 
this is the desired meromorphic continuation of the resolvent 
$(\Delta_h-\lambda)^{-1}$. 
\end{proof}

\noindent
{\bf Examples.}

\noindent 
{\bf 1)} Let $M$ be a surface with cusps. Here by a cusp we mean a 
half-cylinder $[a,\infty)\times S^1$, $a>0$, equipped with the Poincar\'e
metric $y^{-2}(dx^2+dy^2)$, and $M$ is a surface with a complete metric $g$
which in the complement of compact set is isometric to the disjoint union
of finitely many cusps. Let $c>0$ and let $x_0\in M$. Set 
\begin{equation}\label{8.8}
\delta(x):=e^{-cd(x,x_0)}, \quad x\in M,
\end{equation}
and $\rho=\zeta=\delta$. 
Then $\delta$, $\rho$, and $\zeta$ are functions of moderate decay. Let
$$\Omega=\{s\in\C\mid \Re(s)>1/2,\;\; s\notin(1/2,1]\}.$$
We consider
the resolvent $R_g(s)=(\Delta_g-s(1-s))^{-1}$ as a function of $s\in\Omega$.
Then it follows from \cite[Theorem 1]{Mu2} that $R(s)$  admits an
analytic continuation to a meromorphic function on $\C$ with values 
in $\mathcal{L}(L^2_{\delta^{-1}}(M),L^2_\delta(M))$. Using the same method, one
can show that the $R_g(s)$ takes values in $H^2_\delta(M)$. Now observe that the
injectivity radius satisfies $\imath(x)\sim e^{-d(x,x_0)}$. Let $\epsilon>0$
and set $\beta(x)=e^{-(4+\epsilon)d(x,x_0)}$. Choose the constant $c>0$ in
(\ref{8.8}) such that $c<\epsilon/4$. Then $\beta$ is a function of moderate
decay which satisfies (\ref{8.1}) with respect to our choice of the 
functions $\delta$, $\rho$, and
$\zeta$. Now note that the metric $g$ has bounded curvature of all orders.
Let $h$ be complete metric on $M$ with bounded curvature of order 2  
which satisfies $g\sim^2_\beta h$. Then it follows from Theorem \ref{th8.4}
that the resolvent $R_h(s)=(\Delta_h-s(1-s))^{-1}$, $s\in\Omega$, also admits
a meromorphic extension to $\C$ with values in 
$\mathcal{L}(L^2_{\delta^{-1}}(M),H^2_\delta(M))$. We think that the condition
 on $\beta$ can be weakened.  
 
\bigskip
\noindent
{\bf 2)} Let $M$ be a manifold with a cylindrical end. This means that $M$
is a complete Riemannian manifold that admits a decomposition $M=M_0\cup_Y
(\R^+\times Y)$ into a compact manifold $M_0$ with boundary $Y$ and a
half-cylinder $(\R^+\times Y)$ which is glued to $M_0$ along the common 
boundary $Y$. The restriction of the metric $g$ of $M$ to the half-cylinder 
is assumed to be the product metric. Then $g$ is a metric with bounded
geometry, that is, $g$ has bounded curvature of all orders and the injectivity
radius has a positive lower bound. Let $\Delta_Y$ be the Laplacian of $Y$ and
let $0=\mu_1<\mu_2\le \mu_3\le\cdots$ be the eigenvalues of 
$\Delta_Y$. Let $\Sigma\to\C$ be the Riemann surface to which the square 
roots $\lambda\mapsto \sqrt{\lambda-\mu_j}$, $j\in\N$, extend holomorphically.
Define $\delta$, $\rho$, and $\zeta$ as in example 1. Then it follows as in
\cite[Theorem 5]{Mu2} that the resolvent $(\Delta_g-\lambda)^{-1}$
extends from $\C-[0,\infty)$ to a meromorphic function $\lambda\in\Sigma
\mapsto R_g(\lambda)$ with values in 
$\mathcal{L}(L^2_{\delta^{-1}}(M),H^2_\delta(M))$. Now let $\epsilon>0$, 
$x_0\in M$, and set 
$$\beta(x)=e^{-\epsilon d(x,x_0)}, \quad x\in M.$$
Choose 
$c$ in the definition of $\delta$ such that $c<\epsilon/2$. Then $\beta$
satisfies (\ref{8.1}) with respect to our choice of the functions 
$\delta$, $\rho$, and $\zeta$.  Let $h$ be a complete metric on $M$
with bounded curvature of order 2, and suppose that $g\sim^2_\beta h$. Then
it follows from Theorem \ref{th8.4} that the resolvent 
$(\Delta_h-\lambda)^{-1}$
also admits a extension from $\C-[0,\infty)$ to a meromorphic 
function $\lambda\in\Sigma\mapsto R_g(\lambda)$ with values in 
$\mathcal{L}(L^2_{\delta^{-1}}(M),H^2_\delta(M))$.

\end{document}